\documentstyle[amscd,amssymb,array,longtable,url,10pt]{amsart}

\newtheorem{Thm}{Theorem}
\newtheorem{Lem}{Lemma}[section]
\newtheorem{Cor}[Lem]{Corollary}
\newtheorem{Prop}[Lem]{Proposition}

\newtheorem*{ThmI}{Theorem I (Lagarias and Odlyzko \cite{LO1977})}
\newtheorem*{ThmII}{Theorem II (Lagarias, Montgomery, and Odlyzko \cite{LMO1979})}
\newtheorem*{ThmIV}{Theorem III (Lagarias, Montgomery, Odlyzko \cite{LMO1979})}
\newtheorem*{Rk}{Remark}

\begin{document}
\title[]{An explicit upper bound for the least prime ideal in the Chebotarev density theorem}
\author{Jeoung-Hwan Ahn \and Soun-Hi Kwon}
\date{June 13, 2018 }
\maketitle

\begin{center}
\address{Department of Mathematics Education \\
         Korea University \\
         02841, Seoul, Korea \\
         $\{$jh-ahn, sounhikwon$\}$@@korea.ac.kr } \\
\end{center}

\setcounter{equation}{0}
\renewcommand{\theequation}{\arabic{section}.\arabic{equation}}

\begin{abstract}
Lagarias, Montgomery, and Odlyzko proved that there exists an
effectively computable absolute constant $A_1$ such that for every
finite extension $K$ of ${\mathbb{Q}}$, every finite Galois
extension $L$ of $K$ with Galois group $G$ and every conjugacy
class $C$ of $G$, there exists a prime ideal $\mathfrak{p}$ of $K$
which is unramified in $L$, for which
$\left[\frac{L/K}{\mathfrak{p}}\right]=C$, for which
$N_{K/{\mathbb Q}}\,\mathfrak{p}$ is a rational prime, and which
satisfies $N_{K/{\mathbb Q}}\,{\mathfrak{p}} \leq 2 {d_L}^{A_1}$.
In this paper we show without any restriction that $N_{K/{\mathbb
Q}}\,{\mathfrak{p}} \leq {d_L}^{12577}$ if $L \neq {\mathbb Q}$,
using the approach developed by Lagarias, Montgomery, and Odlyzko.
\end{abstract}

\tableofcontents

\section{Introduction}\label{Sec-Introduction}

Let $K$ be a finite algebraic extension of ${\mathbb Q}$, and $L$ a finite Galois extension
of $K$ with Galois group $G$. Let $d_L$ and $d_K$ denote the absolute values of discriminants
of $L$ and $K$, respectively, and let $n_L=[L:\mathbb{Q}]$, $n_K=[K:\mathbb{Q}]$.
To each prime ideal $\mathfrak{p}$ of $K$ unramified in $L$ there corresponds a certain
conjugacy class $C$ of $G$ consisting of the set of Frobenius automorphisms attached to the
prime ideals $\mathfrak{P}$ of $L$ which lie over $\mathfrak{p}$. Denote this conjugacy class
by the Artin symbol $\left[\frac{L/K}{\mathfrak{p}}\right]$. For a conjugacy class $C$ of $G$ let
$$\pi_C (x)=\left| \left\{ {\mathfrak{p}}~|~{\mathfrak{p}}~ {\textrm {a prime ideal of}}~ K,
~{\textrm {unramified in}~L,~\left[\frac{L/K}{\mathfrak{p}}\right]=C,~{\textrm {and}~
N_{K/{\mathbb Q}}\,{\mathfrak{p}} \leq x}} \right\} \right|.$$
The Chebotarev density theorem states that
$$\pi_C (x) \sim \frac{|C|}{|G|} Li(x)$$
as $x \rightarrow \infty$. (See \cite{Hei2010}, \cite{Tsc1926},
\cite{Lang}, \cite{Neu},  and \cite{SL}. See also \cite{Ser1981}
for some extensions of Chebotarev's theorem and applications.) The
error term of this theorem was estimated in \cite{LO1977},
\cite{OEs1979}, and \cite{Win2013}. Lagarias, Montgomery, and
Odlyzko estimated upper bound for the least prime ideal
$\mathfrak{p}$ with $\left[\frac{L/K}{\mathfrak{p}}\right]=C$
under the Generalized Riemann Hypothesis\,(GRH), and
unconditionally, in \cite{LO1977} and \cite{LMO1979},
respectively.

\begin{ThmI}
There exists an effectively computable positive absolute constant $A_0$ such that if the GRH holds for
Dedekind zeta function of $L \neq {\mathbb Q}$, then for every conjugacy class $C$ of $G$ there exists
an unramified prime ideal $\mathfrak{p}$ in $K$ such that $\left[\frac{L/K}{\mathfrak{p}}\right]=C$ and
$$N_{K/{\mathbb Q}}\,{\mathfrak{p}} \leq A_0 (\log d_L)^2.$$
\end{ThmI}

Oesterl\'{e}(\cite{OEs1979}) has stated that if GRH holds, then one may have $A_0=70$.
Bach and Sorenson (\cite{BS1996}) has improved this result in two ways:
If GRH holds, then for any class $C$ of $G$ there is a prime $\mathfrak{p}$ in $K$ of degree 1
over ${\mathbb{Q}}$ with $\left[\frac{L/K}{{\mathfrak{p}}}\right]=C$
and $N_{K/{\mathbb{Q}}} {\mathfrak{p}} \leq (4\log d_L +2.5n_L+5)^2$.
(See also \cite{Bac1990}, \cite{KMurty1994}, and \cite{KMurty}.) Let
$$P(C)= \left\{ {\mathfrak{p}}~|~{\mathfrak{p}}~ {\textrm {a prime ideal of}}~ K,
~{\textrm {unramified in}}~L,~{\textrm {of degree one over}}~{\mathbb Q},~{\textrm {and}}~
\left[\frac{L/K}{\mathfrak{p}}\right]=C \right\}.$$

\begin{ThmII}
There is an absolute, effectively computable constant $A_1$ such that
for every finite extension $K$ of ${\mathbb Q}$,
every finite Galois extension $L$ of $K$, and every conjugacy class $C$ of $G$,
there exists a prime $\mathfrak{p}$ in $P(C)$ which satisfies
$$N_{K/{\mathbb Q}}\,{\mathfrak{p}} \leq 2 {d_L}^{A_1}.$$
\end{ThmII}

See also \cite{Weiss}. When $K=\mathbb Q$ and $L={\mathbb Q}(
{\text{e}}^{2\pi i/q} )$, the conjucacy classes of $G$ correspond
to the residues classes modulo $q$ and Theorem II gives an upper
bound for the least prime in an arithmetic progression
(\cite{LO1977} and \cite{LMO1979}). In this case Theorem II is
weaker than Linnik's theorem (\cite{Lin1944}, \cite{Lin1944bis},
\cite{Bom}). For the least prime in an arithmetic progression, see
for example \cite{Chen1965} - \cite{CL1991}, \cite{GP} -
\cite{HB}, \cite{Hooley} - \cite{Jutila1977}, \cite{Pan1957} -
\cite{Pomerance}, \cite{Wang1986}, \cite{Wang1991}, and
\cite{Xylouris}. If $K=\mathbb Q$, $L={\mathbb Q} ( \sqrt D )$,
and $ \rho$ is the non identity in $Gal(L/\mathbb Q)$, Theorem II
gives an upper bound for the least quadratic nonresidue module
$D$. For this case no upper bound better than Theorem II is known
(\cite{Vinogradov}, \cite{Burgess}, \cite{LO1977}, \cite{LMO1979},
\cite{Ankeny}, \cite{LLS}, \cite{LLS2017}). In this paper we
compute the constant $A_1$.
\begin{Thm}\label{Thm-WithoutGRH}
For every finite extension $K$ of
${\mathbb Q}$, every finite Galois extension $L(\neq {\mathbb Q})$ of $K$ with
Galois group $G$, and every conjugacy class $C$ of $G$, there
exists a prime ideal $\mathfrak{p}$ in $P(C)$ which satisfies
$$N_{K/{\mathbb Q}}\,{\mathfrak{p}} \leq {d_L}^{A_1}$$
with $A_1=12577$.
\end{Thm}

To compute the constant $A_1$ we follow the method developed by
\cite{LMO1979}. In particular, we express zero-free regions for
Dedekind zeta functions, density of zeros of Dedekind zeta
functions, and Deuring-Heilbronn phenomenon with explicit
constants in Sections 5-7 below. Zaman showed in \cite{Zaman} that
$N_{K/{\mathbb Q}}\,{\mathfrak{p}} \ll {d_L}^{40}$ for
sufficiently large $d_L$. See also \cite{TZ}. Winckler proved $A_1
= 27175010$ without any restriction in \cite{Win2015}.
\newline
\newline
\setcounter{equation}{0}

\section{Outline of Lagarias-Montgomery-Odlyzko's method}\label{Sec-OutlineOfProofs}

Let $\Re z$ and $\Im z$ denote the real part and imaginary one of $z \in {\mathbb C}$, respectively.
We review the procedure for the proof of Theorem II in \cite{LMO1979}. Let $g \in C$ and
$$F_C (s)=-\frac{|C|}{|G|} \sum_{\psi}\overline{\psi}(g)\frac{L^{\prime}}{L}(s,\psi,L/K),$$
where $\psi$ runs over the irreducible characters of $G$ and $L(s,\psi,L/K)$ is the Artin L-function
attached to $\psi$. The main parts of \cite{LMO1979} consist of estimates of inverse Mellin transforms
$$\frac{1}{2\pi i} \int_{2-i \infty}^{2+i \infty} F_C (s) k(s)\, ds$$
where $k(s)$ is a kernel function. The main steps of the proof of Theorem II in \cite{LMO1979} are as follows:

\begin{itemize}
\item[$(i)$]
From the orthogonality relations for the characters $\psi$ it follows that for $\Re s>1$
$$F_C (s)=\sum_{\mathfrak{p}} \sum_{m=1}^{\infty}
\theta ({\mathfrak{p}}^m) (\log N_{K/{\mathbb Q}} {\mathfrak{p}} ) (N_{K/{\mathbb Q}} {\mathfrak{p}} )^{-ms}$$
where for prime ideals $\mathfrak{p}$ of $K$ unramified in $L$
$$\theta ({\mathfrak{p}}^m) =
  \begin{cases}
   1 & \textrm {if $\left[\frac{L/K}{\mathfrak{p}}\right]^m=C$,} \\
   0 & \textrm {otherwise,}
  \end{cases}$$
and $|\theta ({\mathfrak{p}}^m)| \leq 1$ if $\mathfrak{p}$ ramifies in $L$.
So we can separate the ${\mathfrak{p}}^m$ with $\left[\frac{L/K}{\mathfrak{p}}\right]^m=C$ from the others.
(See Section 3 of \cite{LO1977}.)
\item[$(ii)$]
Using a method due to Deuring (\cite{Deu1934} and \cite{Mac1968}) $F_C (s)$ can be written as a
linear combination of logarithmic derivatives of Hecke L-functions instead of Artin L-functions.
Let $H=<g>$ be the cyclic subgroup generated by $g$, $E$ the fixed field of $H$.
Then
\begin{eqnarray}\label{Fc-Hecke}
 F_C (s)=-\frac{|C|}{|G|} \sum_{\chi}\overline{\chi}(g)\frac{L^{\prime}}{L}(s,\chi,E),
\end{eqnarray}
where $\chi$ runs over the irreducible characters of $H$, and
$L(s,\chi,E)$ is a Hecke L-function attached to field $E$ with
$\chi({\mathfrak{p}})=\chi\left(\left[\frac{L/E}{{\mathfrak{p}}}\right]\right)$
for all prime ideals $\mathfrak{p}$ of $E$ unramified in L. (See
Section $4$ of \cite{LO1977}.) So, all the singularities of
$F_C(s)$ appear at the zeros and the pole of $\zeta_L (s)$.
\item[$(iii)$] The kernel functions which weight prime ideals of
small norm very heavily are used. Set
$$k_0(s\,;x,y)=\left(\frac{y^{s-1}-x^{s-1}}{s-1}\right)^2 ~~~\textrm{for}~~~y>x>1,$$
$$k_1(s)=k_0(s;x,x^2) ~~~\textrm{for}~~~x\geq 2,$$ and
$$k_2(s)=k_2(s\,;x)=x^{s^2+s}~~~\textrm{for}~~~x\geq 2.$$
In the case that $\zeta_L (s)$ has a real zero very close to 1 we use the kernel $k_2(s)$.
Otherwise we use the kernel $k_1(s)$. The use of the kernel functions is the main innovation of \cite{LMO1979}.
\item[$(iv)$]
For $u>0$ we denote by $\widehat{k}(u)$ the inverse Mellin transform of the kernel function $k(s)$.
Then, for $\Re s >1$,
\begin{align*}
 I &= \frac{1}{2\pi i} \int_{2-i \infty}^{2+i \infty} F_C (s) k(s)\, ds\\
   &= \sum_{\mathfrak{p}} \sum_{m=1}^{\infty} \theta ({\mathfrak{p}}^m) (\log N_{K/{\mathbb Q}} {\mathfrak{p}} )
      \widehat{k} (N_{K/{\mathbb Q}} {\mathfrak{p}} ^{m}),
\end{align*}
where the outer sum is over all prime ideals of $K$. An upper bound ${\mathcal{E}}(\log d_L)$ for
\begin{eqnarray} \label{SumPC-LowerBound1}
 \left|~ I - \sum_{{\mathfrak{p}} \in P(C)} (\log N_{K/{\mathbb Q}} {\mathfrak{p}} )
 \widehat{k}(N_{K/{\mathbb Q}} {\mathfrak{p}} )~ \right|
 \leq {\mathcal{E}}(\log d_L)
\end{eqnarray}
was estimated in (3.15) and (3.16) of \cite{LMO1979}.
\item[$(v)$]
The integral $I$ is evaluated by contour integration:
$$I=\frac{|C|}{|G|}k(1)-\frac{|C|}{|G|} \sum_{\chi} \overline{\chi}(g) \sum_{\rho_{\chi}} k(\rho_{\chi})
+{\mathcal{O}}\left(\frac{|C|}{|G|} n_L k(0) +\frac{|C|}{|G|} k\left(-\frac{1}{2}\right) \log d_L \right),$$
where $\rho_{\chi}$ runs over the zeros of $L(s,\chi,E)$ in the critical strip. (See \cite[Section 3]{LMO1979}.)
So we get
\begin{eqnarray}\label{I-LowerBound}
 \frac{|G|}{|C|} I \geq k(1) - \sum_{\rho} |k(\rho)|
 -c_6 \left\{n_L k(0) + k\left(-\frac{1}{2}\right) \log d_L \right\},
\end{eqnarray}
where $\rho$ runs over the zeros of $\zeta_L (s)$ in the critical strip and $c_6$ is some constant.
Note that $\zeta_L(s)=\prod_{\chi} L(s,\chi,E)$, where $\chi$ runs over the irreducible characters
of $H=Gal(L/E)$. From (\ref{SumPC-LowerBound1}) and (\ref{I-LowerBound}) it follows that
\begin{align}\label{SumPC-LowerBound2}
 \sum_{{\mathfrak{p}} \in P(C)} (\log N_{K/{\mathbb Q}} {\mathfrak{p}} ) \widehat{k}(N_{K/{\mathbb Q}} {\mathfrak{p}} )
 \geq & \,\frac{|C|}{|G|}k(1) - \frac{|C|}{|G|}\sum_{\rho} |k(\rho)| \nonumber \\
 & \,-c_6 \frac{|C|}{|G|}\left\{n_L k(0) + k\left(-\frac{1}{2}\right) \log d_L \right\}-{\mathcal{E}}(\log d_L).
\end{align}
\item[$(vi)$] The sum  $$k(1) - \sum_{\rho} |k(\rho)|$$ is
estimated from below. To do this we need to know the location and
the density of the zeros of $\zeta_L (s)$. If the possible
exceptional zero exists, say $\beta_0$, then $k(\beta_0)$ is
large. The term $k(1)-|k(\beta_0)|$ must be controlled compared to
$\sum_{\rho \neq \beta_0}|k(\rho)|$. We need an enlarged zero-free
region which makes possible $\sum_{\rho \neq \beta_0} |k(\rho)|$
to be small. The Deuring-Heilbronn phenomenon guarantees that the
other zeros of $\zeta_L (s)$ can not be very close to $1$.
\item[$(vii)$] We choose $x$ of the kernel $k(s)$ in terms of
$d_L$ so that the right side of (\ref{SumPC-LowerBound2}) is
positive.
\end{itemize}

\noindent Then Theorem II follows. In the remaining sections of
this paper we will make explicit numerically the constants
intervening in the zero free regions, the density of zeros, and
Deuring-Heilbronn phenomenon of $\zeta_L (s)$, and ultimately
$A_1$.

\section{Prime ideals in $P(C)$}\label{Sec-PrimeIdeals}

In this section we will estimate from above
$$\left|I-\sum_{{\mathfrak{p}}\in P(C)}(\log N_{K/{\mathbb{Q}}} {\mathfrak{p}})
\widehat{k}(N_{K/{\mathbb{Q}}} {\mathfrak{p}})\right|.$$
We will treat carefully their bounds in Section 3 of \cite{LMO1979}.
We begin by recalling the inverse Mellin transform of the kernel functions. They can be easily
computed. For $x \geq 2$ and $u>0$ we have
$$
\widehat{k_1}(u)
= \frac{1}{2\pi i}\int_{a-i\infty}^{a+i\infty}\left\{\frac{x^{2(s-1)}-x^{s-1}}{s-1}\right\}^2u^{-s}\,ds
=\begin{cases}
 u^{-1} \log \frac{x^4}{u} & \textrm {if $\,x^3 \leq u \leq x^4$,}\\
 u^{-1} \log \frac{u}{x^2} & \textrm {if $\,x^2 \leq u \leq x^3$,}\\
 0 & \textrm {otherwise,} \\
\end{cases}
$$
and
$$\widehat{k_2}(u)=\frac{1}{2\pi i}\int_{a-i\infty}^{a+i\infty}x^{s^2+s}u^{-s}\,ds=
(4\pi \log x)^{-\frac{1}{2}} \exp\left\{-\frac{\left(\log \frac{u}{x}\right)^2}{4\log x}\right\},$$
where $a>-\frac{1}{2}$.

\begin{Lem}\label{Lem-SumOverRamifiedPrimes}
Let ${\sum}^{\mathcal{R}}$ denote summation over the prime ideals ${\mathfrak{p}}$ of $K$ that ramify in $L$.
For $x \geq 2$ we have then
\begin{itemize}
\item[$(i)$]
$${\sum}^{\mathcal{R}} \sum_{m=1}^{\infty} \theta ({\mathfrak{p}}^m) (\log N_{K/{\mathbb Q}} {\mathfrak{p}} )
\widehat{k_1} (N_{K/{\mathbb Q}} {\mathfrak{p}}^{m})
\leq \frac{2\log x}{x^2} \log d_L;$$
\item[$(ii)$]
$${\sum}^{\mathcal{R}} \sum_{m \geq 1 \atop N_{K/{\mathbb Q}} {\mathfrak{p}}^{m} \leq x^{5}}
\theta ({\mathfrak{p}}^m) (\log N_{K/{\mathbb Q}} {\mathfrak{p}} ) \widehat{k_2} (N_{K/{\mathbb Q}} {\mathfrak{p}}^{m})
\leq \frac{5}{2\sqrt{\pi}\log 3}\left(\log x\right)^{\frac{1}{2}} \log d_L.$$
\end{itemize}
\end{Lem}

\begin{pf}
\begin{itemize}
\item[$(i)$]
Let ${\mathfrak{p}}$ be a prime ideal of $K$ that is ramified in $L$.
Note that $N_{K/{\mathbb Q}}{\mathfrak{p}} \geq 2$ and
${\sum}^{\mathcal{R}} \log N_{K/{\mathbb Q}} {\mathfrak{p}} \leq \log d_L$.
We have
\begin{align*}
 {\sum}^{\mathcal{R}} \sum_{m=1}^{\infty} \theta ({\mathfrak{p}}^m) (\log N_{K/{\mathbb Q}} {\mathfrak{p}} )
  \widehat{k_1} (N_{K/{\mathbb Q}} {\mathfrak{p}}^{m})
 &\leq \log x {\sum}^{\mathcal{R}} \log N_{K/{\mathbb Q}} {\mathfrak{p}}
  \sum_{m \geq 1 \atop N_{K/{\mathbb Q}} {\mathfrak{p}}^{m} \geq x^2} (N_{K/{\mathbb Q}} {\mathfrak{p}}^{m})^{-1}\\
 &\leq \log x {\sum}^{\mathcal{R}} \frac{\log N_{K/{\mathbb Q}} {\mathfrak{p}}}
  {N_{K/{\mathbb Q}} {\mathfrak{p}}^{m_{{\mathfrak{p}}}}} \left(\frac{1}{1-N_{K/{\mathbb Q}} \mathfrak{p}^{-1}}\right)\\
 &\leq \frac{2\log x}{x^2} \log d_L,
\end{align*}
where $m_{{\mathfrak{p}}}=\left\lceil \frac{\log (x^2)}{\log
N_{K/{\mathbb Q}} {\mathfrak{p}}} \right\rceil$.
\item[$(ii)$]
Let $N_{\mathcal R}$ be the number prime ideals of $K$ that are
ramified in $L/K$. Note that $d_L \geq 3^{N_{\mathcal R}}$. (See
Chap. III, IV of \cite{Ser1968}).) We have
\begin{align*}
 {\sum}^{\mathcal{R}} \sum_{m \geq 1 \atop N_{K/{\mathbb Q}} {\mathfrak{p}}^{m} \leq x^{5}}
  \theta ({\mathfrak{p}}^m) (\log N_{K/{\mathbb Q}} {\mathfrak{p}} ) \widehat{k_2} (N_{K/{\mathbb Q}}{\mathfrak{p}}^{m})
 &\leq (4\pi \log x)^{-\frac{1}{2}} {\sum}^{\mathcal{R}} \log N_{K/{\mathbb Q}} {\mathfrak{p}}
  \sum_{m \geq 1 \atop N_{K/{\mathbb Q}} {\mathfrak{p}}^{m} \leq x^{5}} 1 \\
 &\leq (4\pi \log x)^{-\frac{1}{2}} {\sum}^{\mathcal{R}} {5\log x}\\
 &\leq \frac{5}{2\sqrt{\pi}\log 3} (\log x)^{\frac{1}{2}} \log d_L.
\end{align*}
\end{itemize}
\end{pf}

\begin{Lem}\label{Lem-RS}
\begin{itemize}
\item[$(i)$]{\bf (Rosser and Schoenfeld \cite{RS1962})}
For $x>1$,
$$\pi (x) < \alpha_0 \frac{x}{\log x}$$ with $\alpha_0=1.25506$,
where $\pi (x)$ is the number of primes $p$ with $p \leq x$.
\item[$(ii)$]
For $x>1$,
$$S(x) \leq \frac{2\alpha_0}{\log 2} \sqrt{x},$$
where $S(x)$ is the number of prime powers $p^h$ with $h \geq 2$ and $p^h \leq x$.
\item[$(iii)$]
For $x \geq 101$
$$\sum_{p ~{\textrm{prime}} \atop p^h \geq x^2, h \geq 2} p^{-h}
\leq \frac{4.02\alpha_0}{x\log x}.$$
\end{itemize}
\end{Lem}

\begin{pf}
\begin{itemize}
\item[$(i)$]
See Corollary 1 of \cite{RS1962}.
\item[$(ii)$]
We have
$$S(x) \leq \pi \left(\sqrt{x}\right) \frac{\log x}{\log 2} \leq \frac{2\alpha_0}{\log 2}\sqrt{x}$$ by
$(i)$.
\item[$(iii)$]
We have
$$\sum_{p ~{\textrm{prime}} \atop p^h \geq x^2, h \geq 2} p^{-h}
= \sum_{p ~{\textrm{prime}}} \frac{p^{-h_p}}{1-p^{-1}},$$
where $h_p=\max\left(\left\lceil \frac{\log (x^2)}{\log p} \right\rceil,2\right)$
for each prime $p$.
We observe that
$$\sum_{p \leq x} \frac{p^{-h_p}}{1-p^{-1}} \leq \frac{2}{x^2} \pi(x) \leq \frac{2\alpha_0}{x \log x}.$$
For $x \geq 101$
$$\sum_{p > x} \frac{p^{-h_p}}{1-p^{-1}} \leq \sum_{p > x} \frac{p^{-2}}{1-p^{-1}}
\leq \frac{x}{x-1} \sum_{p > x} p^{-2} \leq 1.01 \sum_{p > x} p^{-2}.$$
By using the integration by parts and $(i)$ we estimate $\sum_{p > x} p^{-2}$ from above. Namely,
$$\sum_{p > x} p^{-2} \leq \int_{x}^{\infty}\frac{1}{t^2} d \pi(t)
\leq \int_{x}^{\infty}  \frac{2\pi(t)}{t^3} dt
\leq \int_{x}^{\infty} \frac{2\alpha_0}{t^2 \log t} dt
\leq \frac{2 \alpha_0}{\log x} \int_{x}^{\infty} \frac{dt}{t^2}=\frac{2 \alpha_0}{x \log x}.$$
Hence,
$$\sum_{p ~{\textrm{prime}}} \frac{p^{-h_p}}{1-p^{-1}} \leq \frac{4.02 \alpha_0}{x \log x},$$
which yields $(iii)$.
\end{itemize}
\end{pf}

\begin{Lem}\label{Lem-SumOverNoRationalPrimesSmall}
For $y \leq \infty$,
let ${\sum}_{y}^{\mathcal{P}}$ denote summation over those $({\mathfrak{p}},m)$
for which $N_{K/{\mathbb Q}} {\mathfrak{p}}^{m}$ is not a rational prime
and $N_{K/{\mathbb Q}} {\mathfrak{p}}^{m} \leq y$.
Then
\begin{itemize}
\item[$(i)$] for $x \geq 101$
$${\sum}_{\infty}^{\mathcal{P}} \theta ({\mathfrak{p}}^m) (\log N_{K/{\mathbb Q}} {\mathfrak{p}} )
\widehat{k_1} (N_{K/{\mathbb Q}} {\mathfrak{p}}^{m}) \leq 16.08 \alpha_0 n_K \frac{\log x}{x};$$
\item[$(ii)$]
for $x \geq 10^{10}$
$${{\sum}_{x^5}^{\mathcal{P}}}
\theta ({\mathfrak{p}}^m) (\log N_{K/{\mathbb Q}} {\mathfrak{p}} ) \widehat{k_2} (N_{K/{\mathbb Q}} {\mathfrak{p}}^{m})
\leq \alpha_1 n_K x^{\frac{3}{4}} (\log x)^{\frac{3}{2}}$$
with
$$\alpha_1=\frac{\alpha_0}{3\sqrt{\pi} \log 2} \left( \frac{15}{10^{\frac{47}{2}}\log 10}
+ 7 +\frac{37}{10^{\frac{5}{2}}} \right)=2.4234\cdots.$$
\end{itemize}
\end{Lem}

\begin{pf}
\begin{itemize}
\item[$(i)$] Since for a positive integer $q$ there are at most
$n_K$ distinct prime power ideals ${\mathfrak{p}}^m$ with
$N_{K/{\mathbb{Q}}} {\mathfrak{p}}^m=q$, it follows that
\begin{align*}
 {\sum}_{\infty}^{\mathcal{P}} \theta ({\mathfrak{p}}^m) (\log N_{K/{\mathbb Q}} {\mathfrak{p}} )
  \widehat{k_1} (N_{K/{\mathbb Q}} {\mathfrak{p}}^{m})
 &\leq \log x {\sum}_{\infty}^{\mathcal{P}} (\log N_{K/{\mathbb Q}} {\mathfrak{p}})
  (N_{K/{\mathbb Q}} {\mathfrak{p}}^{m})^{-1}\\
 &\leq 4 (\log x)^2 n_K \sum_{p ~{\textrm{prime}} \atop x^2 \leq p^h \leq x^4, h \geq 2} p^{-h}.
\end{align*}
Hence, by Lemma \ref{Lem-RS} point $(iii)$ we obtain $(i)$.
\item[$(ii)$]
We have
\begin{align*}
 {\sum}_{x^5}^{\mathcal{P}} \theta ({\mathfrak{p}}^m) (\log N_{K/{\mathbb Q}} {\mathfrak{p}} )
  \widehat{k_2} (N_{K/{\mathbb Q}} {\mathfrak{p}}^{m})
 &\leq n_K \sum_{p ~{\textrm{prime}} \atop p^2 \leq p^h \leq x^5} (\log p^h) \widehat{k_2}(p^h) \\
 &\leq n_K \int_{4}^{x^5} (\log u) \,\widehat{k_2}(u) d S(u),
\end{align*}
where $S(u)$ is as Lemma \ref{Lem-RS} point $(ii)$.
According to Lemma \ref{Lem-RS} point $(ii)$, we have $$S(u) \leq \frac{2\alpha_0}{\log 2} \sqrt{u}.$$
Hence,
\begin{align*}
 \int_{4}^{x^5} (\log u) \,\widehat{k_2}(u) d S(u)
 &\leq (\log x^5)\widehat{k_2}(x^5) S(x^5) +
 \int_{4}^{x^5} \widehat{k_2}(u) \left( \frac{\log u \log \frac{u}{x}}{2\log x} - 1 \right) S(u) \, \frac{du}{u} \\
 &\leq \frac{5\alpha_0}{\sqrt{\pi} \log 2} x^{-\frac{3}{2}} (\log x)^{\frac{1}{2}}
 + \int_{\log \frac{4}{x}}^{4 \log x} \widehat{k_2}(x e^t) \left\{\frac{(t+\log x)t}{2\log x}\right\} S(xe^t)\,dt \\
 &\leq \frac{\alpha_0}{3\sqrt{\pi} \log 2} \left(\frac{15}{x^{\frac{9}{4}} \log x}  +7 +\frac{37}{x^{\frac{1}{4}}} \right)
 x^{\frac{3}{4}} (\log x)^{\frac{3}{2}}.
\end{align*}
\end{itemize}
\end{pf}

\begin{Lem}\label{Lem-SumOverNoRationalPrimesLarge}
For $x \geq 2$, we have
$$\sum_{{\mathfrak{p}}} \sum_{m \geq 1 \atop N_{K/{\mathbb{Q}}} {\mathfrak{p}}^m > x^{5}}
\theta ({\mathfrak{p}}^m) (\log N_{K/{\mathbb Q}} {\mathfrak{p}} ) \widehat{k_2} (N_{K/{\mathbb Q}} {\mathfrak{p}}^{m})
\leq \alpha_2 n_K x \left(\log x\right)^{\frac{1}{2}}$$
with $\alpha_2=\frac{5}{\sqrt{\pi}}.$
\end{Lem}

\begin{pf}
We have
\begin{align*}
 \sum_{{\mathfrak{p}}} \sum_{m \geq 1 \atop N_{K/{\mathbb{Q}}} {\mathfrak{p}}^m > x^{5}}
 \theta ({\mathfrak{p}}^m) (\log N_{K/{\mathbb Q}} {\mathfrak{p}} ) \widehat{k_2} (N_{K/{\mathbb Q}} {\mathfrak{p}}^{m})
 &\leq n_K \sum_{p ~{\textrm{prime}} \atop p^h > x^5} (\log p^h) \widehat{k_2}(p^h) \\
 &\leq n_K \int_{x^{5}}^{\infty}(\log u) \, \widehat{k_2}(u) d T(u),
\end{align*}
where $T(u)$ is the number of prime powers $p^h$ with $h \geq 1$
and $p^h \leq u$. Since $T(u) \leq u$ for $u>0$, we have
\begin{align*}
 \int_{x^5}^{\infty} (\log u) \,\widehat{k_2}(u) d T(u)
 &\leq \int_{x^5}^{\infty} \widehat{k_2}(u) \left( \frac{\log u \log \frac{u}{x}}{2\log x} - 1 \right) T(u) \, \frac{du}{u} \\
 &\leq \int_{4\log x}^{\infty} \widehat{k_2}(x e^t) \left\{\frac{(t+\log x)t}{2\log x}-1\right\} T(xe^t)\,dt \\
 &\leq \alpha_2 x (\log x)^{\frac{1}{2}}.
\end{align*}
\end{pf}

From Lemmas \ref{Lem-SumOverRamifiedPrimes}, \ref{Lem-SumOverNoRationalPrimesSmall},
and \ref{Lem-SumOverNoRationalPrimesLarge} we deduce
an upper bound for
$$\left|I_j-\sum_{{\mathfrak{p}}\in P(C)}(\log N_{K/{\mathbb{Q}}} {\mathfrak{p}})
\widehat{k_j}(N_{K/{\mathbb{Q}}} {\mathfrak{p}})\right|$$
for $j=1,2$ as follows.

\begin{Prop}\label{Prop-IjBound}
Let $k_j(s)$ be as above. Let
$$I_j=\frac{1}{2\pi i} \int_{2-i \infty}^{2+i\infty} F_C(s) k_j(s)ds.$$
Assume that $L \neq {\mathbb{Q}}$. Then
\begin{itemize}
\item[$(i)$] for $x \geq 101$
\begin{align}\label{I1Bound}
 \left|I_1-\sum_{{\mathfrak{p}}\in P(C)}(\log N_{K/{\mathbb{Q}}} {\mathfrak{p}})
 \widehat{k_1}(N_{K/{\mathbb{Q}}} {\mathfrak{p}})\right|
 &\leq \frac{2\log x }{x^2} \log d_L + 16.08\alpha_0 n_K \frac{\log x}{x} \nonumber \\
 &\leq \alpha_3 \frac{\log x}{x} \log d_L
\end{align}
with
$$\alpha_3=\frac{2}{101}+\frac{32.16\alpha_0}{\log 3}=36.759\cdots;$$
\item[$(ii)$]
for $x \geq 10^{10}$
\begin{align}\label{I2Bound}
 &\left|I_2-\sum_{{\mathfrak{p}}\in P(C) \atop N_{K/{\mathbb{Q}}} {\mathfrak{p}} \leq x^{5}}
 (\log N_{K/{\mathbb{Q}}} {\mathfrak{p}}) \widehat{k_2}(N_{K/{\mathbb{Q}}} {\mathfrak{p}})\right| \nonumber \\
 \leq & \frac{5}{2\sqrt{\pi} \log 3} (\log x)^{\frac{1}{2}} \log d_L
 + \alpha_1 n_K x^{\frac{3}{4}} (\log x)^{\frac{3}{2}} + \alpha_2 n_K x \left(\log x\right)^{\frac{1}{2}}
 \leq \alpha_4 x (\log x)^{\frac{1}{2}} \log d_L
\end{align}
with
$$\alpha_4=\frac{1}{\log 3}\left(\frac{10^{-9}}{ 4 \sqrt{\pi} }
+\frac{\alpha_1 \log 10}{5 \sqrt{10} }+2\alpha_2\right)=5.4567\cdots.$$
\end{itemize}
\end{Prop}

Note that $d_L \geq 3^{n_L/2}$ for $n_L \geq 2$.
It follows from
the Hermite-Minkowski's inequality $d_L >
\frac{\pi}{3}\left(\frac{3\pi}{4}\right)^{n_L-1}$ for $n_L>1$. For
$n_L=2$, $d_L \geq 3$, and for $n_L \geq 3$,
$\frac{\pi}{3}\left(\frac{3\pi}{4}\right)^{n_L-1}
=\frac{4}{9}\left(\frac{3\pi}{4}\right)^{n_L}> 3^{{n_L}/{2}}.$
(See also p. 140 of \cite{Sta1974} and p. 291 of \cite{LMO1979}.)

\setcounter{equation}{0}

\section{The Contour integral}\label{Sec-ContourIntegral}

In this section we will evaluate the integrals $I_1$ and $I_2$ by contour integration.
We will use $L(s,\chi)$ to denote $L(s,\chi,E)$.
Let ${\mathcal{F}}(\chi)$ be the conductor of $\chi$ and
$A(\chi)=d_E N_{E/{\mathbb{Q}}} {\mathcal{F}}(\chi)$.
Let
$$\delta(\chi) =
\left\{
\begin{array}{cl}
 1 & \mbox{if} \ \mbox{$\chi$ is the principal character,}\\
 0 & \mbox{otherwise.}
\end{array}
\right. $$
We recall that for each $\chi$ there exist non-negative integers $a(\chi)$, $b(\chi)$ such that
$$a(\chi)+b(\chi)=[E:{\mathbb Q}]=n_E,$$
and such that if we define
$$\gamma_{\chi}(s)=\left\{\pi^{-\frac{s}{2}}\Gamma\left(\frac{s}{2}\right)\right\}^{a(\chi)}
\left\{\pi^{-\frac{s+1}{2}}\Gamma\left(\frac{s+1}{2}\right)\right\}^{b(\chi)}$$
and $$\xi (s,\chi)=\{s(s-1)\}^{\delta(\chi)}A(\chi)^{s/2}\gamma_{\chi}(s)L(s,\chi),$$
then $\xi(s,\chi)$ satisfies the functional equation
$$\xi(1-s,\overline{\chi})=W(\chi) \xi(s,\chi),$$
where $W(\chi)$ is a certain constant of absolute value $1$.
Furthermore, $\xi(s,\chi)$ is an entire function of order $1$ and does not vanish at $s=0$.
By Hadamard product theorem we have for every $s \in \mathbb{C}$
$$-\frac{L^{\prime}}{L}(s,\chi)=\frac{1}{2}\log A(\chi) +\delta(\chi)\left(\frac{1}{s}+\frac{1}{s-1}\right)
+\frac{\gamma_{\chi}^{\prime}}{\gamma_{\chi}}(s)-{\mathcal{B}}(\chi)-\sum_{\rho_{\chi} \in Z(\chi)}
\left(\frac{1}{s-\rho_{\chi}}+\frac{1}{\rho_{\chi}}\right),$$
where ${\mathcal{B}}(\chi)$ is some constant and $Z(\chi)$ denotes the set of nontrivial zeros of $L(s,\chi)$.
(See \cite{Sta1974} and \cite{LO1977}.)
According to (2.8) of \cite{Odl1975}
$$\Re {\mathcal{B}}(\chi)=-\sum_{\rho_{\chi} \in Z(\chi)} \Re \frac{1}{\rho_{\chi}}.$$
Hence, for every $s\in \mathbb{C}$
\begin{eqnarray}\label{ExplicitFormula}
 \Re \left\{-\frac{L^{\prime}}{L}(s,\chi)\right\}=\frac{1}{2}\log A(\chi)
 + \delta(\chi)\Re \left(\frac{1}{s}+\frac{1}{s-1}\right)
 +\Re \frac{\gamma_{\chi}^{\prime}}{\gamma_{\chi}}(s)-\sum_{\rho_{\chi} \in Z(\chi)} \Re \frac{1}{s-\rho_{\chi}}.
\end{eqnarray}
For $j=1$, $2$ we have
$$I_j = \frac{|C|}{|G|} \sum_{\chi} \overline{\chi}(g) J_j(\chi) {\mbox{ by (\ref{Fc-Hecke}),}}$$
where $$J_j(\chi)=\frac{1}{2\pi i}
\int_{2-i\infty}^{2+i\infty}-\frac{L^\prime}{L}(s,\chi) k_j (s)
ds.$$ Assume that $T \geq 2$ does not equal the ordinate of any of
the zeros of $L(s,\chi)$. Consider
$$J_j(\chi,T)=\frac{1}{2\pi i} \int_{B(T)}-\frac{L^\prime}{L}(s,\chi) k_j (s) ds$$
for $j=1$, $2$,
where $B(T)$ is the positively oriented rectangle with vertices $2-iT$, $2+iT$, $-\frac{1}{2}+iT$,
and $-\frac{1}{2}-iT$.
By Cauchy's theorem
\begin{eqnarray}\label{Jj-Cauchy}
 J_j(\chi,T)= \delta(\chi) k_j(1) - \left\{a(\chi)-\delta(\chi)\right\} k_j(0)
 - \sum_{ \rho_{\chi} \in Z(\chi) \atop |\Im \rho_\chi|<T} k_j (\rho_\chi)
\end{eqnarray}
for $j=1$, $2$.

\begin{Lem}\label{Lem-VerticalIntegral}
Let
$$V_j(\chi) = \frac{1}{2\pi i} \int_{-\frac{1}{2}+i\infty}^{-\frac{1}{2}-i\infty}
-\frac{L^\prime}{L}(s,\chi) k_j (s) ds$$ for $j=1$, $2$.
Then
\begin{itemize}
 \item[$(i)$] for $x \geq 101$
              $$\left|V_1(\chi)\right| \leq k_1\left(-\frac{1}{2}\right)\{\mu_{1}\log A(\chi)+ n_E \nu_{1}\},$$
              where $\mu_1=0.75296\cdots$ and $\nu_1=19.405\cdots$;
 \item[$(ii)$] for $x \geq 10^{10}$
              $$\left|V_2(\chi)\right| \leq k_2\left(-\frac{1}{2}\right)\{\mu_{2}\log A(\chi)+ n_E \nu_{2}\},$$
              where $\mu_2=0.058787\cdots$ and $\nu_2=1.4793\cdots$.
\end{itemize}
\end{Lem}

\begin{pf}
Let $s=-\frac{1}{2}+it$.
By \cite[Lemme 5.1]{Win2013}
$$\left|-\frac{L^\prime}{L}\left(-\frac{1}{2}+it,\chi\right)\right| \leq \log A(\chi)+n_E v(t),$$
where
$$v(t)=\log \left(\sqrt{\frac{1}{4}+t^2}+2\right)+\frac{19683}{812}.$$
Moreover, for $x \geq 101$
$$\left|k_1\left(-\frac{1}{2}+it\right)\right|
\leq \frac{x^{-3}(1+x^{-\frac{3}{2}})^2}{\frac{9}{4}+t^2}
= k_1\left(-\frac{1}{2}\right) \left(\frac{1+x^{-\frac{3}{2}}}{1-x^{-\frac{3}{2}}}\right)^2 \left(\frac{9}{9+4t^2}\right)
 \leq k_1\left(-\frac{1}{2}\right)v_1(t)$$
with
$v_1(t)=\left(\frac{1+101^{-\frac{3}{2}}}{1-101^{-\frac{3}{2}}}\right)^2 \left(\frac{9}{9+4t^2}\right)$
and for $x \geq 10^{10}$
$$\left|k_2\left(-\frac{1}{2}+it\right)\right|=x^{-\frac{1}{4}-t^2}=k_2\left(-\frac{1}{2}\right)x^{-t^2}
\leq k_2\left(-\frac{1}{2}\right)v_2(t)$$
with $v_2(t)=10^{-10 t^2}$. Hence,
$$\left|\frac{1}{2\pi i} \int_{-\frac{1}{2}+iT}^{-\frac{1}{2}-iT}
-\frac{L^\prime}{L}(s,\chi)k_j(s)ds\right|
\leq \frac{1}{\pi} k_j\left(-\frac{1}{2}\right) \int_{0}^{T}\{\log A(\chi)+n_E v(t)\}v_j(t)dt.$$
Set
$$\mu_{j} = \frac{1}{\pi} \int_{0}^{\infty} v_j(t)dt {\mbox{ and }}
\nu_{j} = \frac{1}{\pi} \int_{0}^{\infty} v(t) v_j(t)dt.$$
The result follows.
\end{pf}

On the two segments from $2 \pm iT$ to $-\frac{1}{2} \pm iT$ we
proceed with the same way as Section 6 of \cite{LO1977}. (See
Section 3 of \cite{LMO1979}, Section 5 of \cite{Win2013}, and
\cite{Lan1927}.) Let
$${\mathcal{H}}_j(T)=\frac{1}{2\pi i} \int_{-\frac{1}{2}}^{-\frac{1}{4}}
\left\{\frac{L^{\prime}}{L}(\sigma+iT,\chi) k_j(\sigma+iT)
-\frac{L^{\prime}}{L}(\sigma-iT,\chi) k_j(\sigma-iT)\right \} d\sigma$$
and
$${\mathcal {H}}_j^{*}(T)=\frac{1}{2\pi i} \int_{-\frac{1}{4}}^{2}
\left\{\frac{L^{\prime}}{L}(\sigma+iT,\chi) k_j(\sigma+iT)
-\frac{L^{\prime}}{L}(\sigma-iT,\chi) k_j(\sigma-iT)\right \} d\sigma.$$
Then
$${\mathcal {H}}_j(T)+{\mathcal {H}}_j^{*}(T)=\frac{1}{2\pi i}\left\{\int_{2+iT}^{-\frac{1}{2}+iT} -\frac{L^\prime}{L}(s,\chi)k_j(s)ds+\int_{-\frac{1}{2}-iT}^{2-iT}
-\frac{L^\prime}{L}(s,\chi)k_j(s)ds\right\}.$$

\begin{Lem}\label{Lem-HorizontalIntegralHj}
For $j=1$, $2$ we have
$${\mathcal {H}}_j(T) \ll |k_j(iT)|(\log A(\chi)+n_E \log T).$$
\end{Lem}

\begin{pf}
Let $s=\sigma \pm iT$ with $-\frac{1}{2} \leq \sigma \leq -\frac{1}{4}$. Then
$$\frac{L^{\prime}}{L}(s,\chi) \ll \log A(\chi) +n_E \log T$$
by \cite[Lemma 6.2]{LO1977} and $k_j(s) \ll |k_j(iT)|$.
The result follows.
\end{pf}

\begin{Lem}\label{Lem-HorizontalIntegralHjStar1}
Let $-\frac{1}{4} \leq \sigma \leq 2$. Then, we have
$$\frac{L^{\prime}}{L}(\sigma \pm iT,\chi)
-\sum_{\rho_{\chi} \in Z(\chi) \atop |\Im \rho_{\chi} \mp T| \leq 1}
\frac{1}{\sigma \pm iT-\rho_{\chi}} \ll \log A(\chi)+n_E \log T.$$
\end{Lem}

\begin{pf}
See \cite[Lemma 5.6]{LO1977}. (See also \cite[Lemma 4.8]{Win2013}.)
\end{pf}

Therefore, for $j=1$, $2$
\begin{eqnarray*}
 & &{\mathcal {H}}_j^{*}(T)-\frac{1}{2\pi i} \int_{-\frac{1}{4}}^{2}
 \left\{ k_j(\sigma+iT) \sum_{\rho_{\chi} \in Z(\chi) \atop |\Im \rho_{\chi} - T| \leq 1}
 \frac{1}{\sigma + iT-\rho_{\chi}}-k_j(\sigma-iT) \sum_{\rho_{\chi} \in Z(\chi) \atop |\Im \rho_{\chi} + T| \leq 1}
 \frac{1}{\sigma - iT-\rho_{\chi}}\right\} d \sigma \\
 & &\ll |k_j(iT)|(\log A(\chi)+n_E \log T)
\end{eqnarray*}
since $k_j(\sigma \pm iT) \ll |k_j(iT)|$ for $-\frac{1}{4} \leq \sigma \leq 2$.

\begin{Lem}\label{Lem-HorizontalIntegralHjStar2}
Let $\rho_{\chi} \in Z(\chi)$ with $t \neq \Im \rho_{\chi}$. If $|t| \geq 2$, then
$$\int_{-\frac{1}{4}}^{2} \frac{k_j(\sigma + it)}{\sigma + it-\rho_{\chi}} d \sigma \ll |k_j(it)|$$
for $j=1,2$.
\end{Lem}

\begin{pf}
Suppose first that $\Im \rho_{\chi} > t$. Let $B_t$ be the positive oriented rectangle with vertices
$2+i(t-1)$, $2+it$, $-\frac{1}{4}+it$, and $-\frac{1}{4}+i(t-1)$. By Cauchy's theorem,
$$\int_{B_t} \frac{k_j(s)}{s-\rho_{\chi}} ds = 0$$
for $j=1,2$. However, on the three sides of the rectangle other than the segment from $-\frac{1}{4}+it$ to $2+it$,
the integrand is majorized by
$$\alpha_5 |k_j(it)|$$
for some positive constant $\alpha_5$ depending on $x$, which proves the result for $\Im \rho_{\chi} >t$.
A similar proof for $\Im \rho_{\chi} < t$ uses the rectangle with vertices
$2+it$, $2+i(t+1)$, $-\frac{1}{4}+i(t+1)$, and $-\frac{1}{4}+it$.
\end{pf}

For $j=1$, $2$ we have
\begin{eqnarray*}
 & &\frac{1}{2\pi i} \int_{-\frac{1}{4}}^{2}
 \left\{k_j(\sigma+iT) \sum_{\rho_{\chi} \in Z(\chi) \atop |\Im \rho_{\chi} - T| \leq 1}
 \frac{1}{\sigma + iT-\rho_{\chi}}-k_j(\sigma-iT) \sum_{\rho_{\chi} \in Z(\chi) \atop |\Im \rho_{\chi} + T| \leq 1}
 \frac{1}{\sigma - iT-\rho_{\chi}}\right\} d \sigma \\
 & & \ll |k_j(iT)|\{n_{\chi}(T)+n_{\chi}(-T)\} \\
 & & \ll |k_j(iT)|(\log A(\chi)+n_E \log T) {\mbox{ by \cite[Lemma 5.4]{LO1977},}}
\end{eqnarray*}
where $n_{\chi}(T)$ denotes the number of zeros $\rho_{\chi} \in Z(\chi)$ with $|\Im \rho_{\chi} - T| \leq 1$.
We may then conclude as follows.

\begin{Lem} \label{Lem-HorizontalIntegralHjStar}
For $j=1$, $2$ we have
$${\mathcal {H}}_j^{*}(T) \ll |k_j(iT)|(\log A(\chi)+n_E \log T).$$
\end{Lem}

\begin{Lem}\label{Lem-HorizontalIntegral}
For $j=1$, $2$ we have
$$\lim_{T\rightarrow \infty} \frac{1}{2\pi i}\left\{\int_{2+iT}^{-\frac{1}{2}+iT}
-\frac{L^\prime}{L}(s,\chi)k_j(s)ds+\int_{-\frac{1}{2}-iT}^{2-iT}
-\frac{L^\prime}{L}(s,\chi)k_j(s)ds\right\}=0.$$
\end{Lem}

\begin{pf}
By Lemmas \ref{Lem-HorizontalIntegralHj} and \ref{Lem-HorizontalIntegralHjStar}
$${\mathcal {H}}_j(T)+{\mathcal {H}}_j^{*}(T) \ll |k_j(iT)|\{\log A(\chi)+n_E \log T\}.$$
Since
$$|k_j(iT)| \leq
\begin{cases}
 \frac{9}{4x^2(1+T^2)} & \textrm {if $j=1$,} \\
 x^{-T^2} & \textrm {if $j=2$,}
\end{cases}$$
the result follows.
\end{pf}

Letting $T \rightarrow \infty$ in (\ref{Jj-Cauchy}) and combining this
and Lemmas \ref{Lem-HorizontalIntegral} yield
$$J_j(\chi)+V_j(\chi) = \delta(\chi) k_j(1)-\left\{a(\chi)-\delta(\chi)\right\}k_j(0)
- \sum_{\rho_{\chi} \in Z(\chi)} k_j(\rho_{\chi})$$
for $j=1$, $2$. Hence, we have
\begin{align*}
 \frac{|G|}{|C|} I_j
 &= \sum_{\chi} \overline{\chi}(g)J_j(\chi) \\
 &= k_j(1)-k_j(0) \sum_{\chi} \overline{\chi}(g)\left\{a(\chi)-\delta(\chi)\right\}
  -\sum_{\chi} \overline{\chi}(g)\left(\sum_{\rho_{\chi} \in Z(\chi)} k_j(\rho_{\chi})\right)
  -\sum_{\chi} \overline{\chi}(g)V_j(\chi)
\end{align*}
for $j=1$, $2$. Note that by the conductor-discriminant formula
(Chap. VI, Section 3 of \cite{Ser1968})
$$\sum_{\chi} \log A(\chi)= \log d_L.$$ We therefore conclude as
follows.

\begin{Prop}\label{Prop-IjLowerBound}
For $j=1,2$ we have
\begin{eqnarray}\label{Ij-LowerBound}
 \frac{|G|}{|C|} I_j \geq k_j(1)-\sum_{\rho \in Z\left(\zeta_L\right)} |k_j(\rho)|
 - \mu_{j} k_j\left(-\frac{1}{2}\right)\log d_L-n_L\left\{ k_j(0)
 + \nu_{j} k_j\left(-\frac{1}{2}\right) \right\}
\end{eqnarray}
where $Z\left(\zeta_L\right)$ denotes the set of all nontrivial zeros of $\zeta_L(s)$, $\mu_j$
and $\nu_j$ are as in Lemma \ref{Lem-VerticalIntegral}.
\end{Prop}

\setcounter{equation}{0}

\section{Density of zeros of Dedekind zeta functions}\label{Sec-ZeroDensity}

To begin with, we recall that for every $s \in {\mathbb{C}}$ we
have
\begin{eqnarray}\label{ExplicitFormulaForDedekindZetaFunction}
 \Re \left\{-\frac{\zeta_L^{\prime}}{\zeta_L}(s)\right\}=\frac{1}{2}\log d_L
 + \Re \left(\frac{1}{s}+\frac{1}{s-1}\right) +\Re \frac{\gamma_L^{\prime}}{\gamma_L}(s)
 -\sum_{\rho \in Z\left(\zeta_L\right)} \Re \frac{1}{s-\rho},
\end{eqnarray}
where
$$\gamma_L(s)=\left\{\pi^{-\frac{s}{2}}\Gamma\left(\frac{s}{2}\right)\right\}^{r_1+r_2}
\left\{\pi^{-\frac{s+1}{2}}\Gamma\left(\frac{s+1}{2}\right)\right\}^{r_2},$$
$r_1$ and $2r_2$ are the numbers of real and complex embeddings of
$L$. (See Lemma 5.1 of \cite{LO1977} or \cite{Sta1974}.)

For any real number $t$ we let
$$n_L(t)=|\{ \rho=\beta +i \gamma \,|\,
\zeta_L(\rho)=0 \mbox{ with } 0<\beta<1 \mbox{ and } |\gamma-t| \leq 1 \}|.$$
For any complex number $s$ and positive real number $r>0$ we let
$$n(r;s)=|\{ \rho \in Z\left(\zeta_L\right)\,|\,|\rho-s| \leq r\}|.$$
From (\ref{ExplicitFormula}) Lagarias and Odlyzko deduced that
$$n_\chi(t) \ll \log A(\chi)+n_E \log (|t|+2)$$
for all $t$. (See Lemma 5.4 of \cite{LO1977}.)
In this section we will bound $n_L(t)$ and $n(r;s)$ from above using (\ref{ExplicitFormula}).
To do this we need some lemmas.

\begin{Lem}\label{Lem-SumInExplicitFormula}
Let $s=\sigma+it$ with $\sigma>1$. We have
$$\sum_{\rho \in Z\left(\zeta_L\right)} \Re \frac{1}{s-\rho} \geq f_0(\sigma) n_L(t),$$
where
$$f_0(\sigma)=\frac{1}{2}\min\left\{\frac{\sigma-1}{(\sigma-1)^2+1},
\frac{\sigma-\frac{1}{2}}{\left(\sigma-\frac{1}{2}\right)^2+1}\right\}
+\frac{1}{2}\min\left\{\frac{\sigma-\frac{1}{2}}{\left(\sigma-\frac{1}{2}\right)^2+1},
\frac{\sigma}{\sigma^2+1}\right\}.$$
\end{Lem}

\begin{pf}
We have
$$\sum_{\rho \in Z\left(\zeta_L\right)} \Re \frac{1}{s-\rho}
\geq \frac{1}{2}\sum_{\beta+i\gamma \in Z\left(\zeta_L\right) \atop |t-\gamma|\leq 1}
\left\{\frac{\sigma-\beta}{(\sigma-\beta)^2+1} + \frac{\sigma+\beta-1}{(\sigma+\beta-1)^2+1}\right\}
\geq f_0(\sigma) n_L(t).$$
\end{pf}

\begin{Lem}\label{Lem-PolyZetaL}
If $\Re s=\sigma>1$, then
$$\Re \frac{\zeta_L^{\prime}}{\zeta_L}(s) \leq n_L f_1(\sigma),$$
where
$$f_1(\sigma)=-\frac{\zeta_{\mathbb{Q}}^{\prime}}{\zeta_{\mathbb{Q}}}(\sigma).$$
\end{Lem}

\begin{pf}
For $\Re s>1$,
$$-\frac{\zeta_L^{\prime}}{\zeta_L}(s)
=\sum_{\mathfrak{P}}\frac{\log N {\mathfrak{P}}}{{ N {\mathfrak{P}}}^s-1}
=\sum_{{\mathfrak{P}}} \log N {\mathfrak{P}} \sum_{m=1}^{\infty} \, N {\mathfrak{P}}^{-ms},$$
where $\mathfrak{P}$ runs over all prime ideals of $L$.
Comparing $-\frac{\zeta_L^{\prime}}{\zeta_L}(\sigma)$
with $-\frac{\zeta_{\mathbb{Q}}^{\prime}}{\zeta_{\mathbb{Q}}}(\sigma)$ yields
$$\Re \frac{\zeta_L^{\prime}}{\zeta_L}(s) \leq \left|-\frac{\zeta_L^{\prime}}{\zeta_L}(s)\right|
\leq -\frac{\zeta_L^{\prime}}{\zeta_L}(\sigma)
\leq n_L\left\{-\frac{\zeta_{\mathbb{Q}}^{\prime}}{\zeta_{\mathbb{Q}}}(\sigma)\right\}.$$
(See Lemma 3.2 of \cite{LO1977}.)
\end{pf}

See also \cite{Delange}, Lemma $(a)$ of
\cite{Lou1992}, Lemma 3.2 of \cite{Win2013}, p.184 in \cite{EM},
and Proposition 2 of \cite{LouAppear}.

\begin{Lem}\label{Lem-PolyGamma}
Assume that $\Re s > \frac{1}{2}$. We have
\begin{itemize}
\item[$(i)$]
$$\Re \frac{\Gamma^{\prime}}{\Gamma}(s) \leq \log |s| + \frac{1}{3} \leq \alpha_6 \log (|s|+2)$$
with $ \alpha_6=1.08$;
\item[$(ii)$]
$$\Re \frac{\Gamma^{\prime}}{\Gamma}(s) \geq \log |s| - \frac{4}{3} \geq \log (|s|+2)-\alpha_7$$
with $ \alpha_7=\frac{4}{3}+\log 5=2.9427\cdots$.
 \end{itemize}
\end{Lem}

\begin{pf}
For $\Re s >0$,
$$\frac{\Gamma^{\prime}}{\Gamma}(s)=\log s - \frac{1}{2s}
-2\int_{0}^{\infty}\frac{\upsilon}{(s^2+\upsilon^2)(e^{2\pi \upsilon}-1)}\,d\upsilon.$$
(See p. 251 of \cite{WW1996}.) Since $|s^2+\upsilon^2| \geq (\Re s)^2$, we have
$$\left|\int_{0}^{\infty}\frac{\upsilon}{(s^2+\upsilon^2)(e^{2\pi \upsilon}-1)}\,d\upsilon\right|
\leq \frac{1}{(\Re s)^2} \int_{0}^{\infty}\frac{\upsilon}{e^{2\pi \upsilon}-1}\,d\upsilon=\frac{1}{24(\Re s)^2}.$$
If $\Re s > \frac{1}{2}$, then
$$\Re \frac{\Gamma^{\prime}}{\Gamma}(s) \leq \log |s| + \frac{1}{12} \frac{1}{(\Re s)^2}
\leq \log |s| + \frac{1}{3}$$
and
$$\Re \frac{\Gamma^{\prime}}{\Gamma}(s) \geq \log |s| - \frac{1}{2|s|} - \frac{1}{12} \frac{1}{(\Re s)^2}
\geq \log |s| - \frac{4}{3}.$$
Set $\varphi_1(\upsilon)=\alpha_6 \log (\upsilon+2) - \log \upsilon - \frac{1}{3}$ for $\upsilon > \frac{1}{2}$.
Then,
$$\varphi_1^{\prime}(\upsilon)=\frac{(\alpha_6-1)\upsilon-2}{\upsilon(\upsilon+2)}\, \mbox{ and } \varphi_1(\upsilon)>\varphi_1\left(\frac{2}{\alpha_6-1}\right)>0.$$
Hence
$$\Re \frac{\Gamma^{\prime}}{\Gamma}(s) \leq \alpha_6 \log (|s|+2).$$
Set $\varphi_2(\upsilon)=\log \upsilon - \frac{4}{3} - \log (\upsilon+2) + \alpha_7$
for $\upsilon > \frac{1}{2}$.
Then
$$\varphi_2^{\prime}(\upsilon)>0\, \mbox{ and } \varphi_2(\upsilon) > \varphi_2\left(\frac{1}{2}\right)=0.$$
Hence
$$\Re \frac{\Gamma^{\prime}}{\Gamma}(s) \geq \log (|s|+2)-\alpha_7.$$
\end{pf}

\begin{Lem}\label{Lem-PolyGammaL}
Let $s=\sigma+it$. If $\sigma > 1$, then
$$\Re \frac{\gamma_L^{\prime}}{\gamma_L}(s)
\leq n_L \left\{f_2(\sigma) \log(|t|+2)-\frac{1}{2}\log \pi \right\},$$
where
$$f_2(\sigma)=\frac{\alpha_6}{2}\left\{\frac{\log (\sigma+5)}{\log 2}-1\right\}.$$
\end{Lem}

\begin{pf}
By definition and $(i)$ of Lemma \ref{Lem-PolyGamma} we have
\begin{align*}
 \Re \frac{\gamma_L^{\prime}}{\gamma_L}(s)
 &= \frac{(r_1+r_2)}{2} \Re \frac{\Gamma^{\prime}}{\Gamma}\left(\frac{s}{2}\right) + \frac{r_2}{2}
  \Re \frac{\Gamma^{\prime}}{\Gamma}\left(\frac{s+1}{2}\right)-\frac{n_L}{2}\log \pi \\
 &\leq \alpha_6\frac{(r_1+r_2)}{2}  \log \left(\frac{|s|}{2}+2\right)
  +\alpha_6\frac{r_2}{2}  \log \left(\frac{|s+1|}{2}+2\right)-\frac{n_L}{2}\log \pi \\
 &\leq \frac{n_L}{2}\left\{\alpha_6\log \left(\frac{|s+1|}{2}+2\right)-\log \pi \right\}.
\end{align*}
It is sufficient to verify that
\begin{eqnarray}\label{Ineq-PolyGammaL}
 \log \left(\frac{|s+1|}{2}+2\right) \leq \left(\frac{\log(\sigma+5)}{\log 2}-1\right) \log (|t|+2).
\end{eqnarray}
Note that $|s+1| \geq 2|t|$ if and only if $|t| \leq (\sigma+1)/\sqrt{3}$.
If $|t| \geq (\sigma+1)/\sqrt{3}$, then (\ref{Ineq-PolyGammaL}) holds.
We suppose now that $|t| < (\sigma+1)/\sqrt{3}$.
Set $\varphi_3(\upsilon)=\varphi_5(\upsilon)/\varphi_4(\upsilon)$ with $\varphi_4(\upsilon)=\upsilon+2$ and $\varphi_5(\upsilon)=2+\sqrt{(\sigma+1)^2+\upsilon^2}/2$.
Then $\varphi_3^{\prime}(\upsilon) \leq 0$
and $\varphi_5(\upsilon) \leq \left(\frac{\varphi_5(0)}{\varphi_4(0)}\right) \varphi_4(\upsilon)$
for $0 \leq \upsilon < (\sigma+1)/\sqrt{3}$.
For $0 \leq \upsilon < (\sigma+1)/\sqrt{3}$ we have then
$$\frac{\log \varphi_5(\upsilon)}{\log \varphi_4(\upsilon)}
\leq \frac{\log \varphi_4(\upsilon)+\log \varphi_5(0)-\log \varphi_4(0)}{\log \varphi_4(\upsilon)}
\leq \frac{\log \varphi_5(0)}{\log \varphi_4(0)}=\frac{\log(\sigma+5)}{\log 2}-1,$$
which yields (\ref{Ineq-PolyGammaL}).
\end{pf}

We are now ready to bound $n_L(t)$.

\begin{Prop}\label{Prop-ZeroDensityTypeOne}
For all $t$ we have
\begin{eqnarray}\label{nLt-WithoutGRHLong}
 n_L(t) \leq 1.1 \log d_L + 2.09 \log \left\{(|t|+2)^{n_L}\right\} + 0.56 n_L + 4.05.
\end{eqnarray}
In particular, if $L \neq {\mathbb{Q}}$, then
\begin{eqnarray}\label{nLt-WithoutGRHShort}
 n_L(t) \leq 2.72 \log \left\{d_L (|t|+2)^{n_L}\right\}.
\end{eqnarray}
\end{Prop}

\begin{pf}
Combining (\ref{ExplicitFormula}), Lemmas \ref{Lem-SumInExplicitFormula}, \ref{Lem-PolyZetaL}, \ref{Lem-PolyGamma},
and \ref{Lem-PolyGammaL} yields
$$f_0(\sigma)n_L(t) \leq \frac{1}{2}\log d_L + \frac{1}{\sigma}+\frac{1}{\sigma-1}
+n_L\left\{f_2(\sigma)\log (|t|+2) -\frac{1}{2} \log \pi + f_1(\sigma)\right\}$$
for $\sigma>1$. We write
\begin{eqnarray}\label{nLt-UpperBound}
 n_L(t) \leq a_1(\sigma) \log d_L + a_2(\sigma) \log \left\{(|t|+2)^{n_L}\right\} + a_3(\sigma) n_L + a_4(\sigma)
\end{eqnarray}
for $\sigma>1$, where
$$a_1(\sigma)=\frac{1}{2f_0(\sigma)}, ~a_2(\sigma)=\frac{f_2(\sigma)}{f_0(\sigma)},
~a_3(\sigma)=\frac{1}{f_0(\sigma)}\left\{f_1(\sigma)-\frac{1}{2} \log \pi \right\},$$
and
$$a_4(\sigma)=\frac{1}{f_0(\sigma)}\left(\frac{1}{\sigma}+\frac{1}{\sigma-1}\right).$$
We choose now appropriate $\sigma$. If $\sigma=(3+\sqrt{17})/4$,
then (\ref{nLt-UpperBound}) yields (\ref{nLt-WithoutGRHLong}).
For the proof of (\ref{nLt-WithoutGRHShort}), we choose $\sigma=2.45$.
In this case, $a_3(\sigma)<0$  and $2 a_3(\sigma)  + a_4(\sigma)>0$.
Since $n_L \geq 2$, it follows from (\ref{nLt-UpperBound}) that
\begin{align*}
 n_L(t)
 &\leq a_1(\sigma) \log d_L + a_2(\sigma) \log \left\{(|t|+2)^{n_L}\right\} + 2 a_3(\sigma)  + a_4(\sigma)\\
 &\leq B_1 \log d_L + B_2 \log \left\{(|t|+2)^{n_L}\right\},
\end{align*}
where $B_1=a_1(\sigma)+\frac{1}{\log 3}\left\{2 a_3(\sigma)  + a_4(\sigma)\right\}=2.6885\cdots$
and $B_2=a_2(\sigma)=2.7106\cdots$. So, we obtain (\ref{nLt-WithoutGRHShort}).
\end{pf}

See also \cite{KN2012}, \cite{Tru2015}, and Lemme 4.6 of \cite{Win2013}.

\begin{Prop}\label{Prop-ZeroDensityTypeTwo}
Let $r$ be a positive real number.
\begin{itemize}
\item[$(i)$]
Assume that
$$n_L(t) \leq \alpha_8 \log \left\{d_L(|t|+2)^{n_L}\right\}$$
for some $\alpha_8>0$. Then we have
$$n(r;\sigma+it) \leq \alpha_8(1+r)\log \left\{d_L(|t|+r+2)^{n_L}\right\}.$$
\item[$(ii)$]
Assume that $L \neq {\mathbb{Q}}$. If $\sigma \geq 1$ and $0<r\leq 1$, then
$$n(r;\sigma+it) \leq 10\left[1+\frac{2f_2(2)}{5} r \log \left\{d_L(|t|+2)^{n_L}\right\} \right].$$
\end{itemize}
\end{Prop}

\begin{pf} Set
$$Z(r;s)=\{\rho \in Z\left(\zeta_L\right)\,|\,|\rho-s| \leq r\}
~~\textrm {~and~}~~Z(t)=\{\beta+i\gamma \in Z\left(\zeta_L\right)\,|\,|\gamma-t| \leq 1\}.$$
Note that $n(r;s)=|Z(r;s)|$ and $n_L(t)=|Z(t)|$.
\begin{itemize}
\item[$(i)$]
Let $t_1,t_2,\cdots,t_{1+[r]}$ be real numbers such that
$t-r \leq t_1 < \cdots <t_{1+[r]} \leq t+r$ and
$$Z(r;s) \subseteq \bigcup_{i=1}^{1+[r]} Z(t_i).$$
Then
\begin{align*}
 n(r;\sigma+it)
 \leq \sum_{i=1}^{1+[r]} n_L(t_i)
 &\leq \alpha_8 \sum_{i=1}^{1+[r]} \{\log d_L +n_L \log(|t_i|+2)\}\\
 &\leq \alpha_8 (1+r)  \{\log d_L +n_L \log(|t|+r+2)\}.
\end{align*}
\item[$(ii)$]
Write $z=1+r+it$. By (\ref{ExplicitFormula}),
$$\sum_{\rho \in Z\left(\zeta_L\right)} \Re \frac{1}{z-\rho}
= \frac{1}{2} \log d_L +\Re \frac{\gamma_L^\prime}{\gamma_L}(z)+\Re \frac{\zeta_L^\prime}{\zeta_L}(z)
+\Re \left(\frac{1}{z}+\frac{1}{z-1}\right).$$
We now estimate $\Re \frac{\gamma_L^{\prime}}{\gamma_L}(z)$
and $\Re \frac{\zeta_L^{\prime}}{\zeta_L}(z)$ from above.
By Lemma \ref{Lem-PolyGammaL}
$$\Re \frac{\gamma_L^{\prime}}{\gamma_L}(z)
 \leq n_L \left\{f_2(1+r) \log(|t|+2)-\frac{1}{2}\log \pi \right\}
 \leq f_2(1+r) \log \left\{(|t|+2)^{n_L}\right\}.$$
It follows from \cite[Proposition 2]{LouAppear} that
$$\Re \frac{\zeta_L^{\prime}}{\zeta_L}(z)
\leq  \left|\frac{\zeta_L^{\prime}}{\zeta_L}(z)\right|
\leq -\frac{\zeta_L^{\prime}}{\zeta_L}(1+r)
\leq \left(\frac{1-\frac{1}{\sqrt{5}}}{2}\right) \log d_L + \frac{1}{r}.$$
Therefore,
$$\sum_{\rho \in Z\left(\zeta_L\right)} \Re \frac{1}{z-\rho}
\leq \left(1-\frac{1}{2\sqrt{5}}\right)\log d_L +  f_2(1+r) \log \left\{(|t|+2)^{n_L}\right\}
+\frac{2}{r}+\frac{1}{1+r}.$$
Moreover,
$$\sum_{\rho \in Z\left(\zeta_L\right)} \Re \frac{1}{z-\rho}
\geq \sum_{\rho \in Z(2r;z)} \Re \frac{1}{z-\rho}
\geq \frac{1}{4r} n(2r;z).$$
Since $Z(r;\sigma+it) \subseteq Z(r;1+it) \subseteq Z(2r;z)$
and $1-\frac{1}{2\sqrt{5}}<f_2(2)$,
we have
\begin{align*}
 n(r;\sigma+it)
 &\leq n(2r;z) \\
 &\leq 4r\left[\left(1-\frac{1}{2\sqrt{5}}\right)\log d_L
  + f_2(1+r) \log \left\{(|t|+2)^{n_L}\right\} +\frac{2}{r}+\frac{1}{1+r}\right]\\
 &\leq 10 \left[1+\frac{2f_2(2)}{5} r\log \left\{d_L(|t|+2)^{n_L}\right\} \right].
\end{align*}
\end{itemize}
\end{pf}

\setcounter{equation}{0}

\section{Zero-free regions for Dedekind zeta functions}\label{Sec-ZeroFreeRegion}

We abbreviate $N_{L/{\mathbb Q}}$ to $N$.
The classical argument to obtain a zero-free region for $\zeta_L(s)$ starts from (\ref{ExplicitFormula})
and for $\sigma>1$
$$\Re \left[\sum_{m=0}^{d}b_m \left\{-\frac{\zeta_L^{\prime}}{\zeta_L}(\sigma+imt)\right\}\right]
=\Re \sum_{m=0}^{d}b_m \sum_{{\mathfrak{a}}} \frac{\wedge ({\mathfrak{a}})}{N{\mathfrak{a}}^{\sigma+imt}}
\geq 0 $$
where $b_m \geq 0$, ${\mathcal{Q}}(\phi)=\sum_{m=0}^{d}b_m \cos (m\phi) \geq 0$,
$\wedge({\mathfrak{a}})$ is the generalized Von Mangoldt function,
and ${\mathfrak{a}}$ runs over all nonzero ideals of $L$.

Using Stechkin's work one can reduce the constant $\frac{1}{2}$ of
the term $\frac{1}{2} \log A(\chi)$ in (\ref{ExplicitFormula}) to
$\frac{1}{2}\left(1-\frac{1}{\sqrt{5}}\right)$, which yields
larger zero-free regions for $\zeta_L(s)$. (See \cite{Ste1970},
\cite{RS1975}, \cite{Graham}, \cite{MCu1984}, \cite{HB},
\cite{Kad05}, \cite{Kad2012}, \cite{MT2015}, \cite{Lou2015Int},
\cite{Lou2015Ram}, \cite{LouAppear}, and \cite{AK2014}.) It is
known that if $L \neq {\mathbb{Q}}$, then $\zeta_L(s)$ has at most
one zero $\rho=\beta+i\gamma$ with
\begin{eqnarray}\label{ZeroFreeRegion-AK2014}
 \beta>1-\frac{1}{2\log d_L} \mbox{ and } |\gamma|<\frac{1}{2\log d_L}.
\end{eqnarray}
If this zero exists then it must be real and simple.
(See Lemma 3 of \cite{Sta1974}, Lemma 2 of \cite{Hof1980}, and \cite{AK2014}.)
This possible zero is called the exceptional zero and denoted by $\rho_0$.
In this section we will show the following:

\begin{Prop}\label{Prop-ZeroFreeRegion}
Assume that $L \neq {\mathbb{Q}}$.
Let $\rho=\beta+i\gamma$ be a nontrivial zero of $\zeta_L(s)$
with $\rho \neq \rho_0$ and $\tau=|\gamma|+2$.
Then
\begin{eqnarray}\label{ZeroFreeRegion}
 1-\beta>\left(29.57 \log d_L \tau^{n_L}\right)^{-1}.
\end{eqnarray}
\end{Prop}

For the zero-free regions of $\zeta_L(s)$ see also Theorem 1.1 of
\cite{Kad2012}, Lemme 7.1 of \cite{Win2013}, and \cite{Zaman2016}.

We use the Stechkin's work (\cite{Ste1970}) as \cite{MCu1984} and
\cite{Kad2012} and use the same notations as \cite{MCu1984} and
\cite{Kad2012}. Set
$$s=\sigma+it,~\sigma_1=\frac{1+\sqrt{1+4 \sigma^2}}{2},~s_1
=\sigma_1+it,~\kappa=\frac{1}{\sqrt{5}},$$
and
$${\mathbb{F}}(s,z)=\Re \left\{\frac{1}{s-z}+\frac{1}{s-(1-\overline{z})}\right\}.$$
For $\sigma>1$
$$\Re \left\{-\frac{\zeta_L^{\prime}}{\zeta_L}(s)+\kappa \frac{\zeta_L^{\prime}}{\zeta_L}(s_1)\right\}
=\sum_{{\mathfrak{a}}} \frac{\wedge({\mathfrak{a}})}{N
{\mathfrak{a}}^\sigma} \left(1-\frac{\kappa}{N
{\mathfrak{a}}^{\sigma_1-\sigma}}\right) \Re \left(N
{\mathfrak{a}}^{-it}\right),$$ where ${\mathfrak{a}}$ runs over
all nonzero ideals of $L$. Moreover, by (\ref{ExplicitFormula})
\begin{align*}
 \Re \left\{-\frac{\zeta_L^{\prime}}{\zeta_L}(s)+\kappa \frac{\zeta_L^{\prime}}{\zeta_L}(s_1)\right\}
 = & \,\frac{1-\kappa}{2} \log d_L + \Re \left\{\frac{\gamma_L^{\prime}}{\gamma_L}(s)-\kappa \frac{\gamma_L^{\prime}}{\gamma_L}(s_1)\right\}+\{{\mathbb{F}}(s,1)-\kappa {\mathbb{F}}(s_1,1)\}\\
 & \, - {\sum_{\Re \rho \geq \frac{1}{2}}}^{\prime} \{{\mathbb{F}}(s,\rho)-\kappa
 {\mathbb{F}}(s_1,\rho)\},
\end{align*}
where
$${\sum_{\Re \rho \geq \frac{1}{2}}}^{\prime}=\frac{1}{2}
\sum_{\rho \in Z\left(\zeta_L\right) \atop \Re \rho =\frac{1}{2}}
+\sum_{\rho \in Z\left(\zeta_L\right) \atop \frac{1}{2}< \Re \rho \leq 1}.$$

Assume that $b_m \geq 0$ and ${\mathcal{Q}}(\phi)=\sum_{m=0}^{d} b_m \cos (m\phi) \geq 0$.
Then, for $\sigma > 1$
$$\sum_{m=0}^{d}b_m \Re \left\{ -\frac{\zeta_L^{\prime}}{\zeta_L}(\sigma+im\gamma)
+\kappa \frac{\zeta_L^{\prime}}{\zeta_L}(\sigma_1+im\gamma)\right\}
=\sum_{{\mathfrak{a}}} \frac{\wedge ({\mathfrak{a}})}{N {\mathfrak{a}}^{\sigma}}
\left(1-\frac{\kappa}{N {\mathfrak{a}}^{\sigma_1-\sigma}}\right)
{\mathcal{Q}} (\gamma \log N {\mathfrak{a}})\geq 0.$$
So,
\begin{eqnarray}\label{SumSi-Positivity}
 0 \leq S_2 + S_3(\sigma,\gamma)+S_4(\sigma,\gamma)-S_1(\sigma,\gamma),
\end{eqnarray}
where
\begin{eqnarray}\label{S1}
 S_1(\sigma,\gamma)=\sum_{m=0}^{d} b_m {\sum_{\Re \rho \geq \frac{1}{2}}}^{\prime}
 \{{\mathbb{F}}(\sigma+im\gamma,\rho)-\kappa {\mathbb{F}}(\sigma_1+im\gamma,\rho)\},
\end{eqnarray}
\begin{eqnarray}\label{S2}
 S_2=\frac{1-\kappa}{2} {\mathcal{Q}}(0) \log d_L,
\end{eqnarray}
\begin{eqnarray}\label{S3}
 S_3(\sigma,\gamma)=\sum_{m=0}^{d} b_m \{{\mathbb{F}}(\sigma+im\gamma,1)
 -\kappa {\mathbb{F}}(\sigma_1+im\gamma,1)\},
\end{eqnarray}
and
\begin{eqnarray}\label{S4}
 S_4(\sigma,\gamma)=\sum_{m=0}^{d} b_m \Re \left\{ \frac{\gamma_L^{\prime}}{\gamma_L}(\sigma+im\gamma)
 -\kappa \frac{\gamma_L^{\prime}}{\gamma_L}(\sigma_1+im\gamma)\right\}.
\end{eqnarray}

Our proof of Proposition \ref{Prop-ZeroFreeRegion} consists of three parts:
We estimate $S_1(\sigma,\gamma)$ from below, $S_3(\sigma,\gamma)$ and $S_4(\sigma,\gamma)$ from above.
Note that if $\rho$ is a nontrivial zero with $|\gamma| < (2\log d_L)^{-1}$,
then (\ref{ZeroFreeRegion}) is satisfied.
So, we may assume that $\rho \in Z\left(\zeta_L\right)$
and $|\gamma| \geq (2\log d_L)^{-1}$.
Assume that
$$1-\beta \leq (b \log d_L \tau^{n_L})^{-1},$$
where $b \geq 4$ is a constant that will be specified later.
Let $\epsilon=(b\log 12)^{-1}$ and $\sigma-1=(b \log d_L \tau^{n_L})^{-1}.$
That is, $1-\beta \leq \epsilon$ and $\sigma-1 \leq \epsilon$ with $\epsilon \leq (4\log12)^{-1}$.

\begin{Lem}\label{Lem-Ste1970}{\bf (Stechkin \cite{Ste1970})} Let $s=\sigma+it$ with $\sigma>1$.
\begin{itemize}
\item[$(i)$]
If $0<\Re z<1$, then
$${\mathbb{F}}(s,z)-\kappa {\mathbb{F}}(s_1,z) \geq 0.$$
\item[$(ii)$]
If $\Im z=t$ and $\frac{1}{2} \leq \Re z<1$, then
$$\Re \frac{1}{s-1+\overline{z}}-\kappa {\mathbb{F}}(s_1,z) \geq 0.$$
\end{itemize}
\end{Lem}

\begin{Lem}\label{Lem-S1}
Keeping the above notation we have
\begin{eqnarray}\label{S1-LowerBound}
 S_1(\sigma,\gamma) \geq \frac{b_1}{\sigma-\beta}-\{{\mathcal{Q}}(0)-b_1\}\alpha_{10}+\sum_{m \neq 1}
 \frac{b_m(\sigma-\beta)}{(\sigma-\beta)^2+\{(m-1)\gamma\}^2}
\end{eqnarray}
where
$$\alpha_9=\frac{\sqrt{5}-1}{2} {\mbox{ and }}
\alpha_{10}=\kappa \left\{\frac{2 \epsilon}{\alpha_9^2}+\frac{\epsilon}{(\alpha_9^{-1}-\epsilon)^2}\right\}
+\frac{\epsilon}{(1-\epsilon)^2}.$$
\end{Lem}

\begin{pf}
By Lemma \ref{Lem-Ste1970} $(i)$
\begin{eqnarray} \label{S1-LowerBound1}
 S_1(\sigma,\gamma) \geq \sum_{m=0}^{d} b_m \{{\mathbb{F}}(\sigma+im\gamma,\beta+i\gamma)
 -\kappa {\mathbb{F}}(\sigma_1+im\gamma,\beta+i\gamma)\}.
\end{eqnarray}
When $m=1$, we have
\begin{eqnarray}\label{S1-LowerBound2}
 {\mathbb{F}}(\sigma+i\gamma,\beta+i\gamma)-\kappa {\mathbb{F}}(\sigma_1+i\gamma,\beta+i\gamma)
 \geq \frac{1}{\sigma-\beta}
\end{eqnarray}
by Lemma \ref{Lem-Ste1970} $(ii)$. When $m \neq 1$, we have
\begin{eqnarray}\label{S1-LowerBound3}
 & &{\mathbb{F}}(\sigma+im\gamma,\beta+i\gamma)
  -\kappa {\mathbb{F}}(\sigma_1+im\gamma,\beta+i\gamma) \nonumber \\
 & &=\frac{\sigma-\beta}{(\sigma-\beta)^2+\{(m-1)\gamma\}^2}
 -{\mathcal{G}}(\sigma_1-\beta,\sigma_1-1+\beta,\sigma-1+\beta;(m-1)\gamma),
\end{eqnarray}
where
$${\mathcal{G}}(\omega_1,\omega_2,\omega_3;\upsilon)
=\kappa\left(\frac{\omega_1}{\omega_1^2+\upsilon^2}+\frac{\omega_2}{\omega_2^2+\upsilon^2}\right)
-\frac{\omega_3}{\omega_3^2+\upsilon^2}.$$
Note that
\begin{eqnarray}\label{S1-LowerBound4}
 0<\sigma_1-\beta-\alpha_9 \leq 2 \epsilon,~~-\epsilon \leq \sigma_1-1+\beta-\alpha_9^{-1}
 \leq \epsilon, \mbox{ and } -\epsilon \leq \sigma-1+\beta-1 \leq \epsilon.
\end{eqnarray}
For $u>0$ and $u_0>0$
\begin{eqnarray}\label{S1-LowerBound5}
 \left|\frac{u}{u^2+\upsilon^2}-\frac{u_0}{u_0^2+\upsilon^2}\right|
 \leq \frac{|u-u_0|}{\min (u,u_0)^2}.
\end{eqnarray}
(See the proof of Lemma 2.2 of \cite{Kad2012} or that of Lemma 5 of \cite{KN2012}.)
Using (\ref{S1-LowerBound4}), (\ref{S1-LowerBound5}), and the fact that
${\mathcal{G}}(\alpha_9,\alpha_9^{-1},1;\upsilon) \leq 0$
for all $\upsilon \in {\mathbb{R}}$ (\cite[Lemma 2.2 point $(i)$]{Kad2012}
or \cite[Lemma 5 point $(i)$]{KN2012})
we get
\begin{eqnarray}\label{S1-LowerBound6}
 {\mathcal{G}}(\sigma_1-\beta,\sigma_1-1+\beta,\sigma-1+\beta;(m-1)\gamma) \leq \alpha_{10}.
\end{eqnarray}
Substituting (\ref{S1-LowerBound2}), (\ref{S1-LowerBound3}), and (\ref{S1-LowerBound6})
into (\ref{S1-LowerBound1}) yields (\ref{S1-LowerBound}).
\end{pf}

\begin{Lem}\label{Lem-S3}
Keeping the above notation we have
\begin{eqnarray}\label{S3-UpperBound}
 S_3(\sigma,\gamma) \leq \frac{b_0}{\sigma-1}+b_0 f_3(1+\epsilon)
 -\{{\mathcal{Q}}(0)-b_0\}({\mathcal{G}}_0 - \alpha_{11}) +\sum_{m \neq
 0}\frac{b_m(\sigma-1)}{(\sigma-1)^2+(m\gamma)^2},
\end{eqnarray}
where
$$f_3(\sigma)=\frac{1}{\sigma}-\kappa \left(\frac{1}{\sigma_1-1}+\frac{1}{\sigma_1}\right), ~\alpha_{11}=\kappa\left(\frac{\epsilon}{\alpha_9^2}+\frac{\epsilon}{\alpha_9^{-2}}\right)+\epsilon
=(3\kappa+1)\epsilon,$$ and ${\mathcal{G}}_0=-0.121585107$.
\end{Lem}

\begin{pf}
When $m=0$, we have
\begin{eqnarray}\label{S3-UpperBound1}
 {\mathbb{F}}(\sigma,1)-\kappa {\mathbb{F}}(\sigma_1,1)
 = \frac{1}{\sigma-1}+f_3(\sigma) \leq \frac{1}{\sigma-1}+f_3(1+\epsilon)
\end{eqnarray}
since $f_3(\sigma)$ is increasing for $1<\sigma<1.75$.
When $m \neq 0$, we have
\begin{eqnarray}\label{S3-UpperBound2}
 {\mathbb{F}}(\sigma+im\gamma,1)-\kappa {\mathbb{F}}(\sigma_1+im\gamma,1)
 = \frac{\sigma-1}{(\sigma-1)^2+(m\gamma)^2}-{\mathcal{G}}(\sigma_1-1,\sigma_1,\sigma;m\gamma).
\end{eqnarray}
Note that
$0<\sigma_1-1-\alpha_9=\sigma_1-\alpha_9^{-1} \leq \epsilon$ and $0<\sigma-1 \leq \epsilon$.
Using Lemma 2.2 of \cite{Kad2012} we get
\begin{eqnarray}\label{S3-UpperBound3}
 {\mathcal{G}}(\sigma_1-1,\sigma_1,\sigma;m\gamma) \geq {\mathcal{G}}_0 - \alpha_{11}.
\end{eqnarray}
On feeding (\ref{S3-UpperBound1}), (\ref{S3-UpperBound2}), and (\ref{S3-UpperBound3}) into (\ref{S3})
we get (\ref{S3-UpperBound}).
\end{pf}
Let
$$D(m) = \left\{
\begin{array}{ccc}
 \frac{1}{4}\{\Gamma_1(1+\epsilon)+\Gamma_0(1+\epsilon)\}-
 \frac{1-\kappa}{2} \log \pi & \mbox{if} & m=0,\\
 f_4(1+\epsilon) \log m +\alpha_{12} & \mbox{if} & m \neq 0, \\
\end{array} \right. $$
where
$$\Gamma_a(s)=\frac{\Gamma^{\prime}}{\Gamma}\left(\frac{s+a}{2}\right)-
\kappa \frac{\Gamma^{\prime}}{\Gamma}\left(\frac{s_1+a}{2}\right),~~f_4(\sigma)
=\frac{\alpha_6-\kappa}{2} \left\{\frac{\log (\sigma+5)}{\log 2}-1\right\},$$
and
$$\alpha_{12}=\frac{\kappa \alpha_7-(1-\kappa) \log \pi}{2}=0.34162\cdots.$$

\begin{Lem}\label{Lem-S4}
Keeping the above notation we have
$$S_4(\sigma,\gamma) \leq \alpha_{13} \log \tau^{n_L}+ \alpha_{14} n_L,$$
where $\alpha_{13}=\{{\mathcal{Q}}(0)-b_0\}f_4(1+\epsilon)$
and $\alpha_{14}=\sum_{m=0}^{d}b_m D(m)$.
\end{Lem}

\begin{pf}
Since $\Gamma_0(\upsilon)$ and $\Gamma_1(\upsilon)$ are monotonically increasing and
$\Gamma_1(\upsilon)>\Gamma_0(\upsilon)$ for $1<\upsilon<2$,
\begin{align*}
 \Re \left\{\frac{\gamma_L^{\prime}}{\gamma_L}(\sigma)
  -\kappa \frac{\gamma_L^{\prime}}{\gamma_L}(\sigma_1)\right\}
 &= \frac{n_L}{2}\Gamma_0(\sigma)+\frac{r_2}{2}\{\Gamma_1(\sigma)
  -\Gamma_0(\sigma)\}-\frac{1-\kappa}{2} n_L \log \pi\\
 &\leq n_L\left\{\frac{1}{4}\Gamma_1(\sigma)+\frac{1}{4}\Gamma_0(\sigma)
  -\frac{1-\kappa}{2} \log \pi\right\} \leq n_L D(0).
\end{align*}
Set $s=\sigma+im\gamma$ and $s_1=\sigma_1+im\gamma$.
For $m \geq 1$
\begin{align*}
 &\Re \left\{\frac{\gamma_L^{\prime}}{\gamma_L}(s)-\kappa \frac{\gamma_L^{\prime}}{\gamma_L}(s_1)\right\}
 \leq \frac{(r_1+r_2)}{2}\left\{\alpha_6\log \left(\frac{|s|}{2}+2\right)
  -\kappa \log \left(\frac{|s_1|}{2}+2\right)+\kappa \alpha_7\right\}\\
 & \,+\frac{r_2}{2}\left\{\alpha_6\log \left(\frac{|s+1|}{2}+2\right)
  -\kappa \log \left(\frac{|s_1+1|}{2}+2\right)+\kappa \alpha_7\right\}
  -\frac{1-\kappa}{2} n_L \log \pi \mbox{ by Lemma } \ref{Lem-PolyGamma}\\
 &\leq \,\frac{n_L}{2}\left\{(\alpha_6-\kappa)\log \left(\frac{|s+1|}{2}+2\right)
  +\kappa \alpha_7-(1-\kappa) \log \pi \right\}\\
 &\leq n_L \{f_4(\sigma) \log (|m\gamma|+2)+\alpha_{12}\} \mbox{ by } (\ref{Ineq-PolyGammaL})\\
 &\leq n_L \{f_4(1+\epsilon) \log (|\gamma|+2) + D(m)\}.
\end{align*}
Hence
$$
 S_4(\sigma,\gamma)
 \leq b_0 n_L D(0) + n_L \sum_{m=1}^{d} b_m \{f_4(1+\epsilon) \log \tau + D(m)\}
 = \alpha_{13} \log \tau^{n_L} + \alpha_{14} n_L.
$$
\end{pf}

Now, Proposition \ref{Prop-ZeroFreeRegion} is ready to be proven.
Combining (\ref{SumSi-Positivity}), (\ref{S2}), Lemmas \ref{Lem-S1}, \ref{Lem-S3}, and \ref{Lem-S4} yields
\begin{align*}
 0 \leq \, &  \frac{1-\kappa}{2} {\mathcal{Q}}(0) \log d_L + \alpha_{13} \log \tau^{n_L}
 + \alpha_{14} n_L + \alpha_{15} + \frac{b_0}{\sigma-1}-\frac{b_1}{\sigma-\beta}
 + \frac{b_1(\sigma-1)}{(\sigma-1)^2+\gamma^2} \\
 & -\frac{b_0(\sigma-\beta)}{(\sigma-\beta)^2+\gamma^2}
 +\sum_{m=2}^{d}b_m \left\{ \frac{(\sigma-1)}{(\sigma-1)^2+(m\gamma)^2}
 - \frac{(\sigma-\beta)}{(\sigma-\beta)^2+\{(m-1)\gamma\}^2}\right\},
\end{align*}
where $\alpha_{15}=b_0 f_3(1+\epsilon)-\{{\mathcal{Q}}(0)-b_0\}({\mathcal{G}}_0-\alpha_{11})
+\{{\mathcal{Q}}(0)-b_1\} \alpha_{10}$.
Since
$$\frac{b_1(\sigma-1)}{(\sigma-1)^2+\gamma^2}-\frac{b_0(\sigma-\beta)}{(\sigma-\beta)^2+\gamma^2}
\leq \frac{(b_1-b_0)(\sigma-1)}{(\sigma-1)^2+\gamma^2}
\leq (b_1-b_0)\left(\frac{4b}{4+b^2}\right) \log d_L$$
and for $m \geq 2$
$$\frac{(\sigma-1)}{(\sigma-1)^2+(m\gamma)^2}
- \frac{(\sigma-\beta)}{(\sigma-\beta)^2+\{(m-1)\gamma\}^2} \leq 0,$$
it follows that
\begin{eqnarray}\label{ZeroFreeRegion-MainInequlityLong}
 0 \leq \alpha_{16} \log d_L + \alpha_{13} \log \tau^{n_L} + \alpha_{14} n_L
 + \alpha_{15} + \frac{b_0}{\sigma-1}-\frac{b_1}{\sigma-\beta}
\end{eqnarray}
with
$$\alpha_{16}=\frac{1-\kappa}{2} {\mathcal{Q}}(0) + (b_1-b_0) \left(\frac{4b}{4+b^2}\right).$$
Let $0 \leq \delta \leq 1$ and $0 \leq \eta \leq 1$.
Note that $d_L \geq 3^{n_L/2}$.
Set
$$B_{11}=\alpha_{16}+ \frac{2\alpha_{14}}{\log 3} \delta+   \frac{\alpha_{15}}{\log 3}\eta,
~~B_{12}=\alpha_{13}+  \frac{\alpha_{14}}{\log 2} (1-\delta)+  \frac{\alpha_{15}}{2\log 2}(1-\eta),$$
and
$$B_{13}= \max\, (B_{11},B_{12}).$$
The inequality (\ref{ZeroFreeRegion-MainInequlityLong}) is replaced by
\begin{eqnarray}\label{ZeroFreeRegion-MainInequlityShort}
 0 \leq B_{13} \log d_L \tau^{n_L} + \frac{b_0}{\sigma-1}-\frac{b_1}{\sigma-\beta}.
\end{eqnarray}
From (\ref{ZeroFreeRegion-MainInequlityShort}) it follows that
$$1- \beta \geq \left(\frac{b_1}{b_0 b+B_{13}}-\frac{1}{b}\right)(\log d_L \tau^{n_L})^{-1}.$$
We choose ${\mathcal{Q}}(\phi)$ with $b_0 < b_1$, $b$, $\delta$, and $\eta$ as follows:
$${\mathcal{Q}}(\phi)=4(1+\cos \phi)(0.51+\cos \phi)^2, b=8.7, \delta=0.66,
\mbox{ and } \eta=0.26,$$
and obtain (\ref{ZeroFreeRegion}).

\setcounter{equation}{0}

\section{The Deuring-Heilbronn phenomenon}\label{Sec-DHPhenomenon}

The Deuring-Heilbronn phenomenon means that if the exceptional
zero of $\zeta_L(s)$ exists then the other zeros of $\zeta_L(s)$
can not be very close to $s=1$. In \cite{LMO1979} Lagarias,
Montgomery, and Odlyzko proved more precisely the following.

\begin{ThmIV}
There are positive, absolute, effectively computable constants $c_7$ and $c_8$
such that if $\zeta_L(s)$ has a real zero $\omega_0 > 0$ then $\zeta_L(\sigma+it) \neq 0$ for
$$\sigma \geq 1- c_8 \frac{\log \left[\frac{c_7}{(1-\omega_0)
\log \left\{d_L (|t|+2)^{n_L}\right\}}\right]}{\log \left\{d_L (|t|+2)^{n_L}\right\}}$$
with the single exception $\sigma+it=\omega_0$.
\end{ThmIV}

See also \cite{Lin1944bis}. In this section we will estimate the
values of $c_7$ and $c_8$ explicitly. We will use a power sum
inequality as \cite{LMO1979}. We begin by recalling the fact that
$(s-1)\zeta_L(s)$ is an entire function of order one. The Hadamard
product theorem says that
$$(s-1)\zeta_L(s)=s^{r_1+r_2-1} e^{a + b s} \prod_{\omega} \left(1-\frac{s}{\omega}\right)e^{s/\omega}$$
for some constants $a$ and $b$, where $\omega$ runs through all the zeros of $\zeta_L(s)$,
$\omega \neq 0$, including the trivial ones, counted with multiplicity.(\cite{Sta1974})
The Euler product for $\zeta_L(s)$ gives
$$-\frac{\zeta_L^{\prime}}{\zeta_L}(s)=\sum_{{\mathfrak{P}}}
\sum_{m=1}^{\infty} (\log N {\mathfrak{P}})\, (N {\mathfrak{P}})^{-ms}$$
for $\Re s >1$, where ${\mathfrak{P}}$ runs over all prime ideals of $L$.
This series is absolutely convergent for $\Re s > 1$.

Suppose that $\zeta_L(s)$ has a real zero $\omega_0>0$.
Differenciating $(2j-1)$ times
the equality
$$
\sum_{{\mathfrak{P}}}\sum_{m=1}^{\infty} (\log N {\mathfrak{P}})\, (N {\mathfrak{P}})^{-ms}=
\frac{1}{s-1}-b - \sum_{\omega}\left(\frac{1}{s-\omega}+\frac{1}{\omega}\right)-\frac{r_1+r_2-1}{s}
$$
yields that for $\Re s > 1$
and $j \geq 1$
\begin{align*}
 &\frac{1}{(2j-1)!}\sum_{{\mathfrak{P}}}
  \sum_{m=1}^{\infty}(\log N {\mathfrak{P}})(\log N {\mathfrak{P}}^m)^{2j-1}(N {\mathfrak{P}})^{-ms}\\
 &=\frac{1}{(s-1)^{2j}}-\frac{1}{(s-\omega_0)^{2j}}
  -\sum_{\omega \in Z\left(\zeta_L\right) \atop \omega \neq \omega_0}
  \frac{1}{(s-\omega)^{2j}}-\sum_{{\textit{\v{m}}}=0}^{\infty}
  \frac{\ell_{\textit{\v{m}}}}{(s+{\textit{\v{m}}})^{2j}},
\end{align*}
where
$$\ell_{\textit{\v{m}}} = \begin{cases}
  r_1+r_2-1 & {\mbox {if ${\textit{\v{m}}}=0$,}} \\
  r_1+r_2 & {\mbox {if  ${\textit{\v{m}}} \neq 0$ is even,}}\\
  r_2 & {\mbox {if ${\textit{\v{m}}}$ is odd.}}
  \end{cases}$$
If $s_0=\sigma_0+it_0$ with $\sigma_0>1$, then
\begin{align}\label{DifferentiateOfPolyZeta}
 &\frac{1}{(2j-1)!}\sum_{{\mathfrak{P}}}
 \sum_{m=1}^{\infty}(\log N {\mathfrak{P}})(\log N {\mathfrak{P}}^m)^{2j-1}
 N {\mathfrak{P}}^{-m\sigma_0}\left\{1+(N {\mathfrak{P}}^m)^{-it_0}\right\}
 +\sum_{{\textit{\v{m}}}=2}^{\infty}\left\{\frac{\ell_{\textit{\v{m}}}}{(\sigma_0+{\textit{\v{m}}})^{2j}}
 +\frac{\ell_{\textit{\v{m}}}}{(s_0+{\textit{\v{m}}})^{2j}}\right\}
 \nonumber \\
 &=\frac{1}{(\sigma_0-1)^{2j}}-\frac{1}{(\sigma_0-\omega_0)^{2j}}
 +\frac{1}{(s_0-1)^{2j}}-\frac{1}{(s_0-\omega_0)^{2j}} -\sum_{n=1}^{\infty}z_n^j,
\end{align}
where $z_n$ is the series of the terms $(\sigma_0-\omega)^{-2}$ and $(s_0-\omega)^{-2}$
for all $\omega \in \{0,-1\} \cup \left(Z(\zeta_L)\backslash \{\omega_0\} \right)$
such that $\omega$ is counted according to its multiplicity and $|z_n|$ is decreasing for $n \geq 1$.
Since the real part of the left side of (\ref{DifferentiateOfPolyZeta}) is nonnegative,
\begin{align}\label{UpperBoundOfSumOfZnPower}
 \Re \sum_{n=1}^{\infty} z_n^{j}
 \leq & \frac{1}{(\sigma_0-1)^{2j}} - \frac{1}{(\sigma_0-\omega_0)^{2j}}
  +\Re \left[\frac{1}{\{(\sigma_0-1)+it_0\}^{2j}}-\frac{1}{\{(\sigma_0-\omega_0)+it_0\}^{2j}}\right].
\end{align}

To evaluate the constants $c_7$ and $c_8$, first, we estimate the
right side of (\ref{UpperBoundOfSumOfZnPower}) from above.

\begin{Lem}\label{Lem-UpperBoundOfftx}
For $\sigma_0 > 1$, $j \geq 1$, and $0 < \upsilon \leq 1$ we let
$$f_5(\sigma_0+it_0,j;\upsilon)
=\Re \left[\frac{1}{\{(\sigma_0-1)+it_0\}^{2j}}-\frac{1}{\{(\sigma_0-\upsilon)+it_0\}^{2j}}\right].$$
Then $$f_5(\sigma_0,j;\omega_0)+f_5(\sigma_0+it_0,j;\omega_0)
\leq \frac{4j(1-\omega_0)}{(\sigma_0-1)^{2j+1}}.$$
\end{Lem}

\begin{pf}
We have
$$f_5(\sigma_0+it_0,j;\upsilon)
=2j \int_{\sigma_0-1}^{\sigma_0-\upsilon} \Re \left\{
\frac{1}{(y+it_0)^{2j+1}}\right\} dy \leq 2j
\frac{1-\upsilon}{(\sigma_0-1)^{2j+1}}.$$ (See (2.43) of
\cite{Win2015}.) The result follows.
\end{pf}

Second, we estimate $\Re \sum_{n=1}^{\infty} z_n^j$ from below
using Theorem 4.2 of \cite{LMO1979}. (See also Theorem 2.3 of
\cite{Zaman}). Set
$${\mathcal{L}}={\mathcal{L}}(s_0)=|z_1|^{-1} \sum_{n=1}^{\infty}|z_n|.$$
According to Theorem
4.2 of \cite{LMO1979} (see also Theorem 2.3 of
\cite{Zaman})
for any $\textit{\v{c}} > 12$,
there exists $j_0$ with $1 \leq j_0 \leq
\textit{\v{c}} {\mathcal{L}}$ such that
\begin{eqnarray}\label{LowerBoundOfSumOfZnPower}
 \Re \sum_{n=1}^{\infty} z_n^{j_0}
 \geq \left(\frac{\textit{\v{c}}-12}{4\textit{\v{c}}}\right)|z_1|^{j_0}.
\end{eqnarray}

Now we estimate $\sum_{n=1}^{\infty}|z_n|$ from above.

\begin{Lem}\label{Lem-UpperBoundOfSumOfZm}
Let $s_0=\sigma_0+it_0$, $z_n$ and $\omega_0$ be as above. Then we have
\begin{eqnarray}\label{UpperBoundOfLntr}
 \sum_{n=1}^{\infty}|z_n|
 \leq B_{17}(\sigma_0) \log d_L + B_{18}(\sigma_0) \log \left\{(|t_0|+2)^{n_L}\right\}
 + B_{19}(\sigma_0) n_L + B_{20}(\sigma_0),
\end{eqnarray}
where $B_{17}(\sigma_0)=2a_1(\sigma_0)$, $B_{18}(\sigma_0)=a_2(\sigma_0)$,
$B_{19}(\sigma_0)=a_2(\sigma_0)\log 2 + 2 a_3(\sigma_0)+\frac{2}{\sigma_0^2}$, and
$B_{20}(\sigma_0)=2 a_4(\sigma_0)-\frac{2}{\sigma_0^2}$
with
$$a_1(\sigma_0)=\frac{1}{2(\sigma_0-1)},~a_2(\sigma_0)=\frac{f_2(\sigma_0)}{\sigma_0-1},
~a_3(\sigma_0)=-\frac{\log \pi}{2(\sigma_0-1)},$$
and
$$a_4(\sigma_0)=\frac{1}{\sigma_0-1}\left(\frac{1}{\sigma_0}+\frac{1}{\sigma_0-1}\right).$$
(Here, $f_2(\sigma_0)$ is as in Section
\ref{Sec-ZeroDensity}.)
\end{Lem}

\begin{pf}
Note that
$$\sum_{n=1}^{\infty}|z_n|
=\sum_{\omega \in Z\left(\zeta_L\right) \atop \omega \neq \omega_0} \frac{1}{|\sigma_0-\omega|^2}
+\sum_{\omega \in Z\left(\zeta_L\right) \atop \omega \neq \omega_0} \frac{1}{|s_0-\omega|^2}
+\frac{\ell_0}{|\sigma_0|^{2}}+\frac{\ell_0}{|s_0|^{2}}
+\frac{\ell_1}{|\sigma_0+1|^{2}}+\frac{\ell_1}{|s_0+1|^{2}}.$$
As $$\frac{\Re s-1}{|s-\omega|^2} \leq \Re \frac{1}{s-\omega}$$
for $s \in {\mathbb{C}}$ and $\omega \in Z\left(\zeta_L\right)$ we have
\begin{align}\label{Sntr-UpperBound}
 \sum_{\omega \in Z\left(\zeta_L\right) \atop \omega \neq \omega_0} \frac{\Re s-1}{|s-\omega|^2}
 &\leq \sum_{\omega \in Z\left(\zeta_L\right)} \Re \frac{1}{s-\omega} \nonumber \\
 &= \frac{1}{2} \log d_L + \Re \left(\frac{1}{s}+\frac{1}{s-1}\right)
 + \Re \frac{\gamma_L^{\prime}}{\gamma_L}(s) + \Re \frac{\zeta_L^{\prime}}{\zeta_L}(s).
\end{align}
Gathering together the bound in Lemma \ref{Lem-PolyGammaL},
the fact that
$\Re \left\{\frac{\zeta_L^{\prime}}{\zeta_L}(\sigma_0)
+\frac{\zeta_L^{\prime}}{\zeta_L}(\sigma_0+it_0)\right\} \leq 0$,
and (\ref{Sntr-UpperBound}) we get
\begin{align*}
\sum_{\omega \in Z\left(\zeta_L\right) \atop \omega \neq \omega_0} \frac{1}{|\sigma_0-\omega|^2}
+\sum_{\omega \in Z\left(\zeta_L\right) \atop \omega \neq \omega_0} \frac{1}{|s_0-\omega|^2}
\leq& 2a_1(\sigma_0) \log d_L + a_2(\sigma_0) \log \left\{(|t_0|+2)^{n_L}\right\} \\
&+\left\{a_2(\sigma_0) \log 2 + 2 a_3(\sigma_0)\right\}n_L+2a_4(\sigma_0).
\end{align*}
Moreover,
$$
\frac{\ell_0}{|\sigma_0|^{2}}+\frac{\ell_0}{|s_0|^{2}}
+\frac{\ell_1}{|\sigma_0+1|^{2}}+\frac{\ell_1}{|s_0+1|^{2}}
\leq \frac{2(r_1+r_2-1)}{\sigma_0^2}+\frac{2r_2}{(\sigma_0+1)^2}
\leq \frac{2}{\sigma_0^2}n_L - \frac{2}{\sigma_0^2}.
$$
The result follows.
\end{pf}

We are now ready to prove the following.

\begin{Thm}\label{Thm-DH}
Suppose that $L \neq {\mathbb{Q}}$ and $\zeta_L(s)$ has a real zero $\omega_0>0$.
Let $\rho=\beta+i\gamma$ be a zero of $\zeta_L(s)$ with $\rho \neq \omega_0$.
\begin{itemize}
\item[$(i)$]
If $L$ is not an imaginary quadratic number field, then
\begin{eqnarray}\label{Ineq-DHi24}
1-\beta \geq c_8 \frac{\log \left\{\frac{c_7}{(1-\omega_0)\log d_L \tau^{n_L}}\right\}}{\log d_L \tau^{n_L}},
\end{eqnarray}
where $\tau=|\gamma|+2$, $c_7=6.7934\cdots \times 10^{-4}$, and
$c_8=16c_7=\frac{1}{92}$. When $L$ is an imaginary quadratic
number field, then (\ref{Ineq-DHi24}) holds with $c_7=5.5803\cdots
\times 10^{-4}$ and $c_8=16c_7=\frac{1}{112}$. \item[$(ii)$] If
$\rho$ is a nontrivial zero of $\zeta_L(s)$, then
(\ref{Ineq-DHi24}) holds with $c_7=8.1168\cdots \times 10^{-4}$
and $c_8=16c_7=\frac{1}{77}$.
\end{itemize}
\end{Thm}

\begin{pf}
\begin{itemize}
\item[$(i)$]
If $L$ is not an imaginary quadratic number field,
then $\zeta_L(s)$ has a zero at $s=0$ and $|z_1|^{-1} \leq \sigma_0^2$.
Setting $t_0=\gamma$ in (\ref{UpperBoundOfLntr}) yields
 $${\mathcal{L}}
 \leq \sigma_0^2\left\{B_{17}(\sigma_0) \log d_L + B_{18}(\sigma_0) \log \tau^{n_L}
 + B_{19}(\sigma_0) n_L + B_{20}(\sigma_0) \right\}.$$
Note that $B_{19}(\sigma_0) \geq 0$ for $\sigma_0 \geq 1.74$.
For $\sigma_0 \geq 1.74$ and $0 \leq \delta,\,\eta \leq 1$, we let
$$B_{22}(\sigma_0,\delta,\eta)=B_{17}(\sigma_0)+ \frac{2B_{19}(\sigma_0)}{\log 3}\delta
+ \frac{B_{20}(\sigma_0)}{\log 3} \eta,$$
$$B_{23}(\sigma_0,\delta,\eta)=B_{18}(\sigma_0)+\frac{B_{19}(\sigma_0)}{\log 2}(1-\delta)
+\frac{B_{20}(\sigma_0)}{2\log 2}(1-\eta),$$ and
$$B_{24}(\sigma_0,\delta,\eta)=\max \{B_{22}(\sigma_0,\delta,\eta),B_{23}(\sigma_0,\delta,\eta)\}.$$
Then we have
$$
 {\mathcal{L}}
 \leq \sigma_0^2 B_{24}(\sigma_0,\delta,\eta)\log d_L\tau^{n_L}
$$
since $d_L \geq 3^{n_L/2}$ and $n_L \geq 2$.
Note that if $\rho \in Z\left(\zeta_L\right)$, then
$|z_1| \geq |\sigma_0+i\gamma-\rho|^{-2} = |\sigma_0-\beta|^{-2}$
and if $\rho \not \in Z\left(\zeta_L\right)$, then $\rho=\beta \leq 0$ and
$|z_1| \geq |\sigma_0|^{-2} \geq |\sigma_0-\beta|^{-2}$.
Thus
$$|z_1|\geq \frac{1}{(\sigma_0-1)^2}\exp\left\{-2\left(\frac{1-\beta}{\sigma_0-1}\right)\right\}$$
and the bound (\ref{LowerBoundOfSumOfZnPower}) yields
$$\Re \sum_{n=1}^{\infty} z_n^{j_0} \geq \left(\frac{\textit{\v{c}}-12}{4\textit{\v{c}}}\right) \frac{1}{(\sigma_0-1)^{2j_0}}\exp\left\{-2j_0\left(\frac{1-\beta}{\sigma_0-1}\right)\right\}.$$
Combining this with (\ref{UpperBoundOfSumOfZnPower}) and the bound in Lemma \ref{Lem-UpperBoundOfftx}
we have
\begin{eqnarray}\label{Ineq-PowerSum}
\left(\frac{\textit{\v{c}}-12}{4\textit{\v{c}}}\right) \frac{1}{(\sigma_0-1)^{2j_0}}
\exp\left\{-2j_0\left(\frac{1-\beta}{\sigma_0-1}\right)\right\}
\leq \frac{4 j_0 (1-\omega_0)}{(\sigma_0-1)^{2j_0+1}}.
\end{eqnarray}
From $j_0 \leq \textit{\v{c}} {\mathcal{L}}
\leq \textit{\v{c}}\sigma_0^2 B_{24}(\sigma_0,\delta,\eta)\log d_L\tau^{n_L}$ it follows that
\begin{eqnarray}\label{Ineq-DHi}
1-\beta \geq c_8(\textit{\v{c}},\sigma_0,\delta,\eta)
\frac{\log \left\{\frac{c_7(\textit{\v{c}},\sigma_0,\delta,\eta)}{(1-\omega_0)\log d_L \tau^{n_L}}\right\}}
{\log d_L \tau^{n_L}},
\end{eqnarray}
where $c_7(\textit{\v{c}},\sigma_0,\delta,\eta)
=\left(\frac{\textit{\v{c}}-12}{8\textit{\v{c}}}\right)c_8(\textit{\v{c}},\sigma_0,\delta,\eta)$
and
$c_8(\textit{\v{c}},\sigma_0,\delta,\eta)
=\frac{\sigma_0-1}{2\textrm{\v{c}}\sigma_0^2 B_{24}(\sigma_0,\delta,\eta)}.$
Choosing $\textit{\v{c}}=24$, $\sigma_0=7.79$, $\delta=1$, and $\eta=1$ we get (\ref{Ineq-DHi24}).
If $L$ is an imaginary quadratic number field, then $\zeta_L(s)$ has a zero at $s=-1$
and $|z_1|^{-1} \leq (\sigma_0+1)^2$.
We have then $${\mathcal{L}} \leq (\sigma_0+1)^2 B_{24}(\sigma_0,\delta,\eta)\log d_L\tau^{n_L}$$
and $j_0 \leq \textit{\v{c}} {\mathcal{L}}
\leq \textit{\v{c}}(\sigma_0+1)^2 B_{24}(\sigma_0,\delta,\eta)\log d_L\tau^{n_L}$.
Moreover,
$$|z_1| \geq |\sigma_0-\beta|^{-2}
\geq \frac{1}{(\sigma_0-1)^2}\exp\left\{-2\left(\frac{1-\beta}{\sigma_0-1}\right)\right\}$$
since $\zeta_L(s)$ does not have a zero at $s=0$.
From (\ref{Ineq-PowerSum}) we get
$$c_8(\textit{\v{c}},\sigma_0,\delta,\eta)
=\frac{\sigma_0-1}{2\textit{\v{c}}(\sigma_0+1)^2 B_{24}(\sigma_0,\delta,\eta)}.$$
Choosing $\textit{\v{c}}=24$, $\sigma_0=12.21$, $\delta=1$, and $\eta=1$ we get the result.
\item[$(ii)$]
We consider $\sum_{n=1}^{\infty}{\widehat{z}_n}^j$ (instead of $\sum_{n=1}^{\infty} z_n^j$
in (\ref{UpperBoundOfSumOfZnPower})),
where $\widehat{z}_n$ is the series of terms $(\sigma_0-\omega)^{-2}$ and $(\sigma_0+it_0-\omega)^{-2}$
for all $\omega \in Z\left( \zeta_L\right) \backslash \{\omega_0\}$
such that $\omega$ is counted according to
its multiplicity and $\left|\widehat{z}_n\right|$ is decreasing for $n \geq 1$.
Since $$\Re \sum_{n=1}^{\infty} {\widehat{z}_n}^j
+\Re \left\{\frac{\ell_0}{\sigma_0^{2j}} + \frac{\ell_0}{(\sigma_0+it_0)^{2j}} + \frac{\ell_1}{(\sigma_0+1)^{2j}}
+ \frac{\ell_1}{(\sigma_0+it_0+1)^{2j}}\right\}=\Re \sum_{n=1}^{\infty} z_n^j$$
and
$$\Re \left\{\frac{1}{(\sigma_0-\omega)^{2j}}
+ \frac{1}{(\sigma_0+it_0-\omega)^{2j}}\right\} \geq 0
~~~\text{for} ~~~~~\omega = 0, -1~~~,$$
\begin{align}\label{UpperBoundOfSumOfZnPower2}
 \Re \sum_{n=1}^{\infty} {\widehat{z}_n}^{j}
 \leq & \frac{1}{(\sigma_0-1)^{2j}} - \frac{1}{(\sigma_0-\omega_0)^{2j}}
  +\Re \left[\frac{1}{\{(\sigma_0-1)+it_0\}^{2j}}-\frac{1}{\{(\sigma_0-\omega_0)+it_0\}^{2j}}\right].
\end{align}
We use the power-sum inequality in Theorem 4.2 of \cite{LMO1979} for $\sum_{n=1}^{\infty}{\widehat{z}_n}^j$.
Set $\widehat{{\mathcal{L}}}=|\widehat{z}_1|^{-1} \sum_{n=1}^{\infty} |{\widehat{z}_n}|$.
For any $\textit{\v{c}} > 12$, there exists $\widehat{j}_0$ with
$1 \leq \widehat{j}_0 \leq {\textit{\v{c}}} \widehat{{\mathcal{L}}}$ such that
\begin{align}\label{LowerBoundOfSumOfZnPower2}
 \Re \sum_{n=1}^{\infty} {\widehat{z}_n}^{\widehat{j}_0}
 \geq \left(\frac{\textit{\v{c}}-12}{4\textit{\v{c}}}\right)|\widehat{z}_1|^{\widehat{j}_0}.
\end{align}
If $\rho \in Z\left(\zeta_L\right)$ , then $1-\overline{\rho} \in
Z\left(\zeta_L\right)$. Set $t_0 = \gamma$. Then
$$|\widehat{z}_1|^{-1}  \leq \min \{(\sigma_0-\beta)^2,
(\sigma_0-1+\beta)^2\} \leq \left(\sigma_0-\frac{1}{2}\right)^2.$$
Then we have
$$\widehat{{\mathcal{L}}} \leq \left(\sigma_0-\frac{1}{2}\right)^2
\left\{B_{17}(\sigma_0)\log d_L + B_{18}(\sigma_0)\log \tau^{n_L}
+ \widehat{B}_{19}(\sigma_0)n_L + \widehat{B}_{20}(\sigma_0)\right\},$$
where $\widehat{B}_{19}(\sigma_0)=a_2(\sigma_0)\log 2 + 2 a_3(\sigma_0)$
and $\widehat{B}_{20}(\sigma_0)=2a_4(\sigma_0)$.
Note that $\widehat{B}_{19}(\sigma_0) \leq 0$
and $2\widehat{B}_{19}(\sigma_0)+\widehat{B}_{20}(\sigma_0) \geq 0$
for $1<\sigma_0 \leq 11.66$. So, for $1 < \sigma_0 \leq 11.66$
$$\widehat{{\mathcal{L}}}
 \leq \left(\sigma_0-\frac{1}{2}\right)^2\left\{B_{17}(\sigma_0)\log d_L
 + B_{18}(\sigma_0)\log \tau^{n_L}
 + 2\widehat{B}_{19}(\sigma_0) + \widehat{B}_{20}(\sigma_0)\right\}.$$
For $1 < \sigma_0 \leq 11.66$ and $0 \leq \eta \leq 1$, we let
$$B_{25}(\sigma_0,\eta)=B_{17}(\sigma_0)+ \frac{2\widehat{B}_{19}(\sigma_0)
+ \widehat{B}_{20}(\sigma_0)}{\log 3} \eta,$$
$$B_{26}(\sigma_0,\eta)=B_{18}(\sigma_0)+\frac{2\widehat{B}_{19}(\sigma_0)
+ \widehat{B}_{20}(\sigma_0)}{2\log 2}(1-\eta),$$
and
$$B_{27}(\sigma_0,\eta)=\max \{B_{25}(\sigma_0,\eta),B_{26}(\sigma_0,\eta)\}.$$
Then we have
$$
 \widehat{{\mathcal{L}}}
\leq \left(\sigma_0-\frac{1}{2}\right)^2 B_{27}(\sigma_0,\eta)\log d_L\tau^{n_L}.
$$
Note that $d_L \geq 3^{n_L/2}$.
Since $|z_1| \geq |\sigma_0+i\gamma-\rho|^{-2}
\geq \frac{1}{(\sigma_0-1)^2}\exp\left\{-2\left(\frac{1-\beta}{\sigma_0-1}\right)\right\}$,
the bound (\ref{LowerBoundOfSumOfZnPower2}) yields
$$\Re \sum_{n=1}^{\infty} \widehat{z}_n^{\widehat{j}_0}
\geq \left(\frac{\textit{\v{c}}-12}{4\textit{\v{c}}}\right) \frac{1}{(\sigma_0-1)^{2\widehat{j}_0}}
\exp\left\{-2\widehat{j}_0\left(\frac{1-\beta}{\sigma_0-1}\right)\right\}.$$
Combining this with (\ref{UpperBoundOfSumOfZnPower2})
and the bound in Lemma \ref{Lem-UpperBoundOfftx} we have
$$\left(\frac{\textit{\v{c}}-12}{4\textit{\v{c}}}\right) \frac{1}{(\sigma_0-1)^{2\widehat{j}_0}}
\exp\left\{-2\widehat{j}_0\left(\frac{1-\beta}{\sigma_0-1}\right)\right\}
\leq \frac{4 \widehat{j}_0 (1-\omega_0)}{(\sigma_0-1)^{2\widehat{j}_0+1}}.$$
From $\widehat{j}_0 \leq \textit{\v{c}} {\mathcal{L}}
\leq \textit{\v{c}}\left(\sigma_0-\frac{1}{2}\right)^2 B_{27}(\sigma_0,\eta)\log d_L\tau^{n_L}$
it follows that
$$c_8(\textit{\v{c}},\sigma_0,\eta)
=\frac{\sigma_0-1}{2\textit{\v{c}}\left(\sigma_0-\frac{1}{2}\right)^2 B_{27}(\sigma_0,\eta)}.$$
Choosing $\textit{\v{c}}=24$, $\sigma_0=5.42$, and $\eta=1$ we get the result.
\end{itemize}
\end{pf}

\begin{Rk}\label{Rk-UpperBoundOfS}
To get an upper bound for ${\mathcal{L}}$
the zero-density estimate for the number of
zeros of $\zeta_L(s)$ was used in \cite{LMO1979}:
\begin{align*}
 {\mathcal{L}}
 &\ll (2-\beta)^2 \sum_{\omega} \left(\frac{1}{|2-\omega|^2}+\frac{1}{|2+i\gamma-\omega|}\right)\\
 &\ll \int_{0}^{\infty} \frac{1}{u^2+1} d n(u) + \int_{0}^{\infty} \frac{1}{u^2+1} d n(u+\tau) \\
 &\ll \log d_L \tau^{n_L},
\end{align*}
where $\omega$ runs through all the zeros of $\zeta_L(s)$ including the trivial ones.
(See (5.6) of \cite{LMO1979}.) However we used
$$\sum_{\rho \in Z\left(\zeta_L\right)} \frac{\sigma-1}{|s-\rho|^2}
\leq \sum_{\rho \in Z\left(\zeta_L\right)} \Re \frac{1}{s-\rho}$$
for $\Re s=\sigma>1$ and (\ref{ExplicitFormulaForDedekindZetaFunction}).
(See (\ref{Sntr-UpperBound}) above.)
\end{Rk}

\begin{Cor}\label{Cor-ZFRORZ}
Assume that $L \neq {\mathbb{Q}}$. Then for any real zero
$\omega_0>0$ of $\zeta_L(s)$ we have
\begin{eqnarray}\label{LowerBoundOfdL}
1-\omega_0 \geq d_L^{-c_{10}}
\end{eqnarray}
with $c_{10}=114.72\cdots$.
\end{Cor}

\begin{pf}
When $L$ is not an imaginary quadratic number fields, we let
$\textit{\v{c}}=12.1$, $\sigma_0=7.79$, $\delta=1$, and $\eta=1$.
The inequality (\ref{Ineq-DHi}) yields
\begin{align}\label{Ineq-NTRZ}
1-\beta \geq c_8 \frac{\log c_7 + \log (1-\omega_0)^{-1}
- \log \log d_L \tau^{n_L}}{\log d_L \tau^{n_L}}
\end{align}
for any zero $\beta+i\gamma \neq \omega_0$ of $\zeta_L(s)$,
where $c_7=2.2434\cdots \times 10^{-5}$ and $c_8=2.1716\cdots \times 10^{-2}$.
Set $1-\omega_0=d_L^{-c}$.
Since $\zeta_L(s)$ always has a trivial zero at $s=0$ and $d_L \geq 3^{n_L/2}$,
we have
\begin{align}\label{Ineq-RZ}
 1
 &\geq c_8 \left\{\frac{\log c_7+c \log d_L}{\left(1+\frac{2\log 2}{\log 3}\right) \log d_L}
 -\frac{\log \log d_L 2^{n_L}}{\log d_L 2^{n_L}}\right\} \nonumber \\
 &\geq c_8 \left\{\left(1+\frac{2\log 2}{\log 3}\right)^{-1}\left(\frac{\log c_7}{\log d_L}+c\right)
 -\frac{1}{e}\right\}.
\end{align}
Note that $\frac{\log x}{x} \leq \frac{1}{e}$ for $x > 0$. Then (\ref{Ineq-RZ}) yields
$$
 c \leq \left(\frac{1}{c_8}+\frac{1}{e}\right)\left(1+\frac{2\log 2}{\log 3}\right)
 -\frac{\log c_7}{\log 3} = 114.72\cdots.
$$
When $L$ is an imaginary quadratic number field,
it is known that $\zeta_L(\sigma) \neq 0$ for $\sigma \geq 1- \left(\frac{\pi}{6}\sqrt{d_L}\right)^{-1}$.
(See the proof of Lemma 11 of \cite{Sta1974}.)
The result follows.
\end{pf}

\begin{Rk}
\begin{itemize}
\item[(1)] For the zero-free regions for $\zeta_L(s)$ see also
\cite{Sta1974}. \item[(2)] In \cite{Zaman}, Zaman proved that, for
$d_L$ sufficiently large, $1-\omega_0 \gg d_L^{-21. 3}$.
\end{itemize}
\end{Rk}

\setcounter{equation}{0}

\section{Proof of Theorem \ref{Thm-WithoutGRH}}\label{Sec-ProofOfMainTheorem}

Theorem \ref{Thm-WithoutGRH} is ready to be proven. We will choose
appropriate kernel functions $k(s)$ and estimate
$$k(1)-\sum_{\rho \in Z\left(\zeta_L\right)} |k(\rho)|$$ from below. From now
on we denote by $\beta_0$ the exceptional zero of $\zeta_L(s)$ if
it exists, and $\beta_0=1-(2\log d_L)^{-1}$ otherwise. Our proof
is divided into a sequence of lemmas.

\begin{Lem}\label{Lem-LowerBoundForTrivialZero}
We have
\begin{eqnarray}\label{LowerBoundForTrivialZeroK1}
 k_1(1)-k_1(\beta_0) \geq \frac{9}{10}(\log x)^2 \min\{1,(1-\beta_0)\log x\}
\end{eqnarray}
and
\begin{eqnarray}\label{LowerBoundForTrivialZeroK2}
 k_2(1)-k_2(\beta_0) \geq \frac{9}{10}x^2 \min\{1,(1-\beta_0)\log x\}.
\end{eqnarray}
\end{Lem}

\begin{pf}
We have
\begin{eqnarray*}
 k_1(1)-k_1(\beta_0)&=&(\log x)^2-\left(\frac{x^{(\beta_0-1)}-x^{2(\beta_0-1)}}{1-\beta_0}\right)^2\\
 &=&(\log x)^2 \varphi_6((1-\beta_0)\log x),
\end{eqnarray*}
where
$$\varphi_6(\upsilon)=1-\left(\frac{e^{-\upsilon}-e^{-2\upsilon}}{\upsilon}\right)^2.$$
It is easily verified that
$$\varphi_6(\upsilon) \geq
 \begin{cases}
  \varphi_6(1) \upsilon & \textrm {for $0 < \upsilon \leq 1$,} \\
  \varphi_6(1) & \textrm {for $\upsilon \geq 1$}
 \end{cases}$$
with $\varphi_6(1)=0.94592\cdots$. Hence $\varphi_6(\upsilon) \geq \varphi_6(1) \min\{1,\upsilon\}$,
which yields (\ref{LowerBoundForTrivialZeroK1}). We have
$$k_2(1)-k_2(\beta_0)=x^2(1-x^{(\beta_0-1)(\beta_0+2)}) \geq x^2 \varphi_7((1-\beta_0)\log x),$$
where $\varphi_7(\upsilon)=1-e^{-\frac{5}{2}\upsilon}$. It is easy to see that
$$\varphi_7(\upsilon) \geq
 \begin{cases}
  \varphi_7(1) \upsilon & \textrm {for $0 < \upsilon \leq 1$,} \\
  \varphi_7(1) & \textrm {for $\upsilon \geq 1$}
 \end{cases}$$
with $\varphi_7(1)=0.91791\cdots$. Hence $\varphi_7(\upsilon) \geq \varphi_7(1) \min\{1,\upsilon\}$,
which yields (\ref{LowerBoundForTrivialZeroK2}).
\end{pf}

In the following $c_7$ and $c_8$ are as in Theorem \ref{Thm-DH} point $(ii)$.

\begin{Lem}\label{Lem-LowerBoundForSumOfNontrivialZerosK1}
Suppose that $\beta_0 \leq 1-{c_7}^2 (\log d_L 3^{n_L})^{-2}$.
We use the kernel function $k_1(s)$ and obtain
$$\sum_{\rho \in Z\left(\zeta_L\right) \atop \rho \neq \beta_0} |k_1(\rho)|
\leq c_{13} \log d_L+ c_{14}(\log d_L)^2\{(1-\beta_0) \log d_L\}^{2c_{12}\frac{\log x}{\log d_L}},$$
where $c_{12}=6.8610\cdots \times 10^{-4}$, $c_{13}=124.14\cdots$, and $c_{14}=1.7700\cdots \times 10^8$.
\end{Lem}

\begin{pf}
Write
$$\sum_{\rho \neq \beta_0 \atop \rho \in Z\left(\zeta_L\right)} |k_1(\rho)|
=\sum_{|\rho-1| > 1 } |k_1(\rho)|+\sum_{|\rho-1| \leq 1 } |k_1(\rho)|,$$
where $\sum_{|\rho-1| > 1}$ (resp. $\sum_{|\rho-1| \leq 1}$) denotes that we sum over $\rho=\beta+i\gamma$ such that
$\rho \in Z\left(\zeta_L\right)$ with $\rho \neq \beta_0$ and $|\rho-1| > 1$ (resp. $|\rho-1| \leq 1$).
Since
$$|k_1(\rho)|=\left|\frac{x^{2(\rho-1)}-x^{\rho-1}}{\rho-1}\right|^2 \leq \frac{4x^{-2(1-\beta)}}{|\rho-1|^2},$$
it follows that
\begin{align*}
 \sum_{|\rho-1| > 1} |k_1(\rho)|
 &\leq 4 \int_{1}^{\infty} \frac{1}{r^2} d n(r;1)\\
 &\leq 21.76 \int_{1}^{\infty}\frac{(1+r)\{\log d_L + n_L \log (r+2)\}}{r^3} dr
 {\mbox{ by (\ref{nLt-WithoutGRHShort}) and Proposition \ref{Prop-ZeroDensityTypeTwo} $(i)$}}\\
 &\leq c_{13} \log d_L
\end{align*}
where $c_{13}=21.76\left(\frac{3}{2}+\frac{2+15\log 3}{4\log 3}\right)=124.14\cdots$.
For the sum $\sum_{|\rho-1| \leq 1} |k_1(\rho)|$ we consider two cases separately.
\begin{itemize}
\item[$(i)$]
If an exceptional zero $\beta_0$ exists with
$1-\beta_0 \leq \left(\frac{c_7}{3}\right)^2 (\log d_L)^{-1}$,
then
$$\frac{c_7}{(1-\beta_0)\log d_L \tau^{n_L}} \geq \frac{c_7}{3(1-\beta_0)\log d_L}
\geq \left\{(1-\beta_0)\log d_L\right\}^{-\frac{1}{2}}.$$
Hence, by Theorem \ref{Thm-DH} point $(ii)$
$$1-\beta \geq c_8 \frac{\log \left\{(1-\beta_0)\log d_L\right\}^{-\frac{1}{2}}}{\log d_L \tau^{n_L}}
\geq c_{11} \frac{\log \left\{(1-\beta_0)\log d_L\right\}^{-1}}{\log d_L}$$
with $c_{11}=\frac{c_8}{6}=\frac{1}{462}$.
\item[$(ii)$]
If $1-\beta_0 > \left(\frac{c_7}{3}\right)^2 (\log d_L)^{-1}$, then by (\ref{ZeroFreeRegion})
$$1-\beta \geq (29.57 \log d_L \tau^{n_L})^{-1} \geq (88.71 \log d_L)^{-1}.$$
Set $c_{12}=\left\{177.42\log \left(\frac{3}{c_7}\right)\right\}^{-1}=6.8610\cdots \times 10^{-4}$.
Then $$(88.71)^{-1}=2 c_{12} \log \left(\frac{3}{c_7}\right)
> c_{12} \log \left\{(1-\beta_0)\log d_L\right\}^{-1}$$ and
$$1-\beta > c_{12}  \frac{\log \left\{(1-\beta_0)\log d_L\right\}^{-1}}{\log d_L}.$$
\end{itemize}
As $c_{11} > c_{12}$ we have
$$1-\beta > c_{12}  \frac{\log \left\{(1-\beta_0)\log d_L\right\}^{-1}}{\log d_L}$$
in all cases. Let
$$B= c_{12}  \frac{\log \left\{(1-\beta_0)\log d_L\right\}^{-1}}{\log d_L}.$$
Then
$$|k_1(\rho)| \leq \frac{4x^{2(\beta-1)}}{|\rho-1|^2} \leq \frac{4x^{-2B}}{|\rho-1|^2}.$$
By Proposition (\ref{Prop-ZeroDensityTypeTwo}) point $(ii)$,
\begin{align*}
 \sum_{|\rho-1| \leq 1} |k_1(\rho)|
 &\leq 4 x^{-2B}\int_{B}^{1} \frac{1}{r^2} d n(r;1)\\
 &\leq 4x^{-2B}\left\{n(1;1)+20\int_{B}^{1}\frac{1+\frac{2f_2(2)}{5}
  \left(1+\frac{2 \log 2}{\log 3}\right)r\log d_L}{r^3}dr\right\}
 {\mbox{ by Proposition \ref{Prop-ZeroDensityTypeTwo} $(ii)$}}\\
 &\leq 40x^{-2B}\left\{B^{-2}+\frac{4f_2(2)}{5}\left(1+\frac{2 \log 2}{\log 3}\right)B^{-1} \log d_L
 - \frac{2f_2(2)}{5}\left(1+\frac{2 \log 2}{\log 3}\right) \log d_L \right\}\\
 &\leq c_{14} (\log d_L)^2\left\{(1-\beta_0)\log d_L\right\}^{2c_{12}\frac{\log x}{\log d_L}}
\end{align*}
where
$$c_{14}=\frac{40}{c_{12}\log 2}\left\{\frac{1}{c_{12}\log 2}
+\frac{4f_2(2)}{5}\left(1+\frac{2 \log 2}{\log 3}\right)\right\}
=1.7700\cdots \times 10^8.$$
For the last inequality we used (\ref{ZeroFreeRegion-AK2014}),
which yields
$$B= c_{12}  \frac{\log \left\{(1-\beta_0)\log d_L\right\}^{-1}}{\log d_L}
\geq \frac{c_{12}\log 2}{\log d_L}.$$
\end{pf}

We have therefore
\begin{align}\label{I1-LowerBoundFront}
 k_1(1)-\sum_{\rho \in Z\left(\zeta_L\right)} |k_1(\rho)|
 \geq & \frac{9}{10} (\log x)^2 \min\{1,(1-\beta_0)\log x\} -c_{13} \log d_L \nonumber \\
 & - c_{14} (\log d_L)^2\left\{(1-\beta_0)\log d_L\right\}^{2c_{12}\frac{\log x}{\log d_L}}.
\end{align}
Note that for $x \geq 101$
\begin{align}\label{I1-LowerBoundBack}
 & \mu_{1} k_1\left(-\frac{1}{2}\right)\log d_L
 +n_L\left\{ k_1(0)+ \nu_{1} k_1\left(-\frac{1}{2}\right) \right\} \nonumber\\
 & \leq \left\{\frac{2}{\log 3}\left(x^{-2}-x^{-1}\right)^2
 +\frac{4}{9}\left(\mu_1+\frac{2}{\log 3}\nu_1\right)\left(x^{-3}-x^{-3/2}\right)^2\right\}\log d_L \nonumber \\
 & \leq \left\{\frac{2}{\log 3} x^{-2}
 + \frac{4}{9}\left(\mu_1+\frac{2}{\log 3}\nu_1\right)x^{-3}\right\} \log d_L \nonumber \\
 & \leq c_{15} x^{-2}\log d_L,
\end{align}
where
$$c_{15}=\frac{2}{\log 3}+\frac{4}{909}\left(\mu_1+\frac{2}{\log 3}\nu_1\right)=1.9792\cdots.$$
Gathering together the bounds (\ref{I1Bound}),
(\ref{Ij-LowerBound}), (\ref{I1-LowerBoundFront}), and
(\ref{I1-LowerBoundBack}) we conclude the following:

\begin{Lem}\label{Lem-LowerBoundForSumPCK1}
Suppose that $\beta_0 \leq 1-{c_7}^2 (\log d_L 3^{n_L})^{-2}$. We have then
\begin{align}\label{LowerBoundForSumPCk1}
 &\frac{|G|}{|C|} \sum_{{\mathfrak{p}} \in P(C)} (\log N_{K/{\mathbb Q}} {\mathfrak{p}} )
 \widehat{k_1}(N_{K/{\mathbb Q}} {\mathfrak{p}} )
 \geq \frac{9}{10} (\log x)^2 \min\{1,(1-\beta_0)\log x\}-c_{13} \log d_L  \nonumber \\
 &- c_{14} (\log d_L)^2\left\{(1-\beta_0)\log d_L\right\}^{2c_{12}
 \frac{\log x}{\log d_L}} -c_{15} x^{-2}\log d_L - \alpha_3 \frac{|G|}{|C|} \frac{\log x}{x}\log d_L.
\end{align}
\end{Lem}

\begin{Lem}\label{Lem-PositivityOfSumOfPCK1}
Suppose that $\beta_0 \leq 1-{c_7}^2 (\log d_L 3^{n_L})^{-2}$. For $\log x=c_{16} \log d_L$
with $c_{16}=3144.25$, we have
$$\sum_{{\mathfrak{p}} \in P(C)} (\log N_{K/{\mathbb Q}} {\mathfrak{p}} )
\widehat{k_1}(N_{K/{\mathbb Q}} {\mathfrak{p}} )>0.$$
In particular, there is a prime ${\mathfrak{p}} \in P(C)$ with $N_{K/{\mathbb Q}} {\mathfrak{p}} \leq x^4 = d_L^{4c_{16}}$.
\end{Lem}

\begin{pf}
Let $\log x=c_{16} \log d_L$.
\begin{itemize}
\item[$(i)$]
Suppose that $1 \leq c_{16}(1-\beta_0)\log d_L$.
(\ref{LowerBoundForSumPCk1}) and (\ref{ZeroFreeRegion-AK2014}) yield
$$(\log d_L)^{-2} \frac{|G|}{|C|} \sum_{{\mathfrak{p}} \in P(C)} (\log N_{K/{\mathbb Q}} {\mathfrak{p}} )
\widehat{k_1}(N_{K/{\mathbb Q}} {\mathfrak{p}} )
\geq \left\{\frac{9}{10}c_{16}^2-c_{14}\left(\frac{1}{2}\right)^{2c_{12}c_{16}}\right\}-\epsilon_1,$$
where
$$\epsilon_1 = \frac{c_{13}}{\log d_L}+\frac{c_{15}}{d_L^{2c_{16}} \log d_L}
+\frac{2 \alpha_3 c_{16}\log d_L}{d_L^{c_{16}}\log 3}.$$
(Note that $\frac{|G|}{|C|} \leq |G| = \frac{n_L}{n_K} \leq n_L \leq \frac{2}{\log 3} \log d_L$.)
For $c_{16}=3144.25$, we have
$$\frac{9}{10} c_{16}^2 > c_{14}\left(\frac{1}{2}\right)^{2c_{12}c_{16}}+\epsilon_1.$$
\item[$(ii)$]
Suppose that $1 \geq c_{16}(1-\beta_0)\log d_L$.
Since $1-\beta_0 \geq c_7^2(\log d_L 3^{n_L})^{-2}
\geq \left(\frac{c_7}{3}\right)^2(\log d_L)^{-2}$,
(\ref{LowerBoundForSumPCk1}) and (\ref{ZeroFreeRegion-AK2014}) yield
\begin{align*}
 & \left\{(1-\beta_0)\log d_L\right\}^{-1}(\log d_L)^{-2} \frac{|G|}{|C|}
 \sum_{{\mathfrak{p}} \in P(C)} (\log N_{K/{\mathbb Q}} {\mathfrak{p}} )
 \widehat{k_1}(N_{K/{\mathbb Q}} {\mathfrak{p}} )\\
 &\geq  \frac{9}{10} c_{16}^3 - c_{14}\left\{(1-\beta_0)\log d_L\right\}^{2c_{12}c_{16}-1}
 -\frac{c_{13}}{(1-\beta_0)(\log d_L)^2} -\frac{c_{15}}{d_L^{2c_{16}} (1-\beta_0) (\log d_L)^2}\\
 &- \frac{2 \alpha_3 c_{16}}{d_L^{c_{16}}(1-\beta_0)\log 3}
 \geq  \frac{9}{10} c_{16}^3 - c_{14} \left(\frac{1}{2}\right)^{2c_{12}c_{16}-1}
 - c_{13}\left(\frac{3}{c_7}\right)^2 - \epsilon_2,
\end{align*}
where
$$\epsilon_2=\left(\frac{3}{c_7}\right)^2 \left\{\frac{c_{15}}{d_L^{2c_{16}}}
+ \frac{2\alpha_3 c_{16}}{\log 3} \frac{(\log d_L)^2}{d_L^{c_{16}}} \right\}.$$
For $c_{16} = 1261$, we have
$$\frac{9}{10} c_{16}^3 > c_{14} \left(\frac{1}{2}\right)^{2c_{12}c_{16}-1}
+c_{13}\left(\frac{3}{c_7}\right)^2 + \epsilon_2.$$
\end{itemize}
The result follows.
\end{pf}

\begin{Lem}\label{Lem-LowerBoundForSumPCK2}
Suppose that $1-\beta_0 \leq {c_7}^2 (\log d_L 3^{n_L})^{-2}$.
We have then
\begin{align}\label{LowerBoundForSumPCk2}
 &\frac{|G|}{|C|} \sum_{{\mathfrak{p}} \in P(C) \atop N_{K/{\mathbb{Q}}} {\mathfrak{p}} \leq x^{5}}
 (\log N_{K/{\mathbb Q}} {\mathfrak{p}} ) \widehat{k_2}(N_{K/{\mathbb Q}} {\mathfrak{p}} )
 \geq \frac{9}{10} x^2 \min\{1,(1-\beta_0)\log x\} - c_{20} x \log d_L \nonumber \\
 & - c_{21} x^2 (1-\beta_0)^{2c_{19}\frac{\log x}{\log d_L}} \log d_L - c_{15}^{\prime} \log d_L
 - \alpha_4 \frac{|G|}{|C|} x(\log x)^{\frac{1}{2}}\log d_L,
\end{align}
where $c_{20}=19.16\cdots$, $c_{21}=6.1522\cdots$, $c_{19}=\frac{c_8}{6}=\frac{1}{462}$, and
$c_{15}^{\prime}=1.8291\cdots$.
\end{Lem}

\begin{pf}
For $\rho = \beta+i\gamma \in Z\left(\zeta_L\right)$ with $|\gamma| \leq 1$
we have by Theorem \ref{Thm-DH} point $(ii)$
$$1-\beta \geq c_8 \frac{\log \left\{\frac{c_7}{(1-\beta_0) \log d_L 3^{n_L}}\right\}}{\log d_L 3^{n_L}}
\geq c_{19} \frac{\log (1-\beta_0)^{-1}}{\log d_L}$$
with $c_{19}=\frac{c_8}{6}=\frac{1}{462}$.
Since
$|k_2(\rho)| \leq x^{\beta^2+\beta} \leq x^{2-2(1-\beta)}
 \leq x^2(1-\beta_0)^{2c_{19}\frac{\log x}{\log d_L}}$,
$$\sum_{|\gamma| \leq 1} |k_2(\rho)|
 \leq x^2(1-\beta_0)^{2c_{19}\frac{\log x}{\log d_L}} \sum_{|\gamma| \leq 1} 1
 \leq c_{21} x^2(1-\beta_0)^{2c_{19}\frac{\log x}{\log d_L}} \log d_L
 {\mbox{ by (\ref{nLt-WithoutGRHShort})}}$$
with $c_{21}=2.72\left(1+\frac{2 \log 2}{\log 3}\right)=6.1522\cdots$.
For zeros $\rho=\beta+i\gamma$ with $|\gamma|>1$ and $x \geq 10^{10}$ we have
\begin{align*}
 \sum_{|\gamma|>1} |k_2(\rho)|
 &\leq x^2 \sum_{m=1}^{\infty} \{n_L(2m)+n_L(-2m)\} x^{-(2m-1)^2} \\
 &\leq 5.44 x^2 \sum_{m=1}^{\infty} \{\log d_L + n_L \log (2m+2)\} x^{-(2m-1)^2}
 {\mbox{ by (\ref{nLt-WithoutGRHShort})}} \\
 &\leq c_{20} x \log d_L,
\end{align*}
where
$$c_{20}=5.44 \sum_{m=1}^{\infty} \left\{1 + \frac{2}{\log 3} \log (2m+2)\right\}
10^{-40m^2+40m}=19.16\cdots.$$
It follows that for $x \geq 10^{10}$
\begin{align}\label{I2-LowerBoundFront}
 k_2(1)-\sum_{\rho} |k_2(\rho)|
 \geq& \frac{9}{10} x^2 \min\{1,(1-\beta_0)\log x\} \nonumber \\
 & -c_{21}x^2(1-\beta_0)^{2c_{19}\frac{\log x}{\log d_L}} \log d_L - c_{20} x \log d_L.
\end{align}
Note that for $x \geq 10^{10}$
\begin{align}\label{I2-LowerBoundBack}
 \mu_{2} k_2\left(-\frac{1}{2}\right)\log d_L
 +n_L\left\{ k_2(0)+ \nu_{2} k_2\left(-\frac{1}{2}\right) \right\}
 &\leq \left\{\frac{2}{\log 3}+\left(\mu_2+\frac{2}{\log 3}\nu_2\right)x^{-\frac{1}{4}}\right\}\log d_L \nonumber \\
 &\leq c_{15}^\prime \log d_L,
\end{align}
where
$$ c_{15}^\prime=\frac{2}{\log 3}+\left(\mu_2+\frac{2}{\log 3}\nu_2\right)10^{-\frac{5}{2}}=1.8291\cdots.$$
Combining (\ref{I2Bound}), (\ref{Ij-LowerBound}),
(\ref{I2-LowerBoundFront}), and (\ref{I2-LowerBoundBack}) yields
(\ref{LowerBoundForSumPCk2}).
\end{pf}

\begin{Lem}\label{Lem-PositivityOfSumOfPCK2}
Suppose that $1-\beta_0 \leq {c_7}^2 (\log d_L 3^{n_L})^{-2}$. If $x=d_L^{c_{23}}$
with $c_{23}=179$, then
$$ \sum_{{\mathfrak{p}} \in P(C) \atop N_{K/{\mathbb{Q}}} {\mathfrak{p}}
\leq x^{5}} (\log N_{K/{\mathbb Q}} {\mathfrak{p}} ) \widehat{k_2}(N_{K/{\mathbb Q}} {\mathfrak{p}} )>0.$$
In particular, there is a prime ${\mathfrak{p}} \in P(C)$
with $N_{K/{\mathbb Q}} {\mathfrak{p}} \leq x^{5} = d_L^{5c_{23}}$.
\end{Lem}

\begin{pf}
Let $x=d_L^{c_{23}}$. Then (\ref{LowerBoundForSumPCk2}) becomes
\begin{align*}
 &\frac{|G|}{|C|} \sum_{{\mathfrak{p}} \in P(C) \atop N_{K/{\mathbb{Q}}} {\mathfrak{p}} \leq x^{5}}
 (\log N_{K/{\mathbb Q}} {\mathfrak{p}} ) \widehat{k_2}(N_{K/{\mathbb Q}} {\mathfrak{p}} )
 \geq \frac{9}{10} d_L^{2c_{23}} \min\{1,c_{23}(1-\beta_0)\log d_L\}  \nonumber \\
 & -c_{20} d_L^{c_{23}} \log d_L - c_{21} d_L^{2c_{23}} (1-\beta_0)^{2c_{19}c_{23}} \log d_L - c_{15}^{\prime} \log d_L
 - \frac{2\alpha_4 c_{23}^{\frac{1}{2}}}{\log 3} d_L^{c_{23}}(\log d_L)^{\frac{5}{2}}.
\end{align*}
When $1 \leq c_{23}(1-\beta_0) \log d_L$, we have
\begin{align*}
 d_L^{-2c_{23}}\frac{|G|}{|C|} \sum_{{\mathfrak{p}} \in P(C) \atop N_{K/{\mathbb{Q}}} {\mathfrak{p}} \leq x^{5}}
 (\log N_{K/{\mathbb Q}} {\mathfrak{p}} ) \widehat{k_2}(N_{K/{\mathbb Q}} {\mathfrak{p}} )
 &\geq \frac{9}{10}-c_{21}\{c_7^2(\log d_L)^{-2}\}^{2c_{19}c_{23}} \log d_L -\epsilon_3 \\
 &= \frac{9}{10}-c_{21} c_7^{4c_{19}c_{23}} (\log d_L)^{1-4c_{19}c_{23}} -\epsilon_3,
\end{align*}
where
$$\epsilon_3=c_{20}\frac{\log d_L}{d_L^{c_{23}}}+c_{15}^{\prime}\frac{\log d_L}{d_L^{2c_{23}}}+
\frac{2\alpha_4 c_{23}^{\frac{1}{2}}}{\log 3}\frac{(\log d_L)^{\frac{5}{2}}}{d_L^{c_{23}}}.$$
If $c_{23}=(4c_{19})^{-1}=114.76\cdots$, then
$$\frac{9}{10}>c_{21}c_7 +\epsilon_3.$$
When $1 \geq c_{23}(1-\beta_0) \log d_L$, using Corollary \ref{Cor-ZFRORZ} we have
\begin{align*}
 & d_L^{-2c_{23}}\{(1-\beta_0)\log d_L\}^{-1}\frac{|G|}{|C|}
 \sum_{{\mathfrak{p}} \in P(C) \atop N_{K/{\mathbb{Q}}} {\mathfrak{p}} \leq x^{5}}
  (\log N_{K/{\mathbb Q}} {\mathfrak{p}} ) \widehat{k_2}(N_{K/{\mathbb Q}} {\mathfrak{p}} ) \\
 &\geq \frac{9}{10}c_{23} - \frac{c_{20}}{d_L^{c_{23}}(1-\beta_0)} - c_{21}(1-\beta_0)^{2c_{19}c_{23}-1}
 - \frac{c_{15}^{\prime}}{d_L^{2c_{23}}(1-\beta_0)} - \frac{2\alpha_4 c_{23}^{\frac{1}{2}}}{\log 3}
 \frac{(\log d_L)^{\frac{3}{2}}}{d_L^{c_{23}}(1-\beta_0)} \\
 & \geq \frac{9}{10}c_{23} - \epsilon_4,
\end{align*}
where
$$\epsilon_4=\frac{c_{20}}{d_L^{c_{23}-c_{10}}} + c_{21}c_7^{4c_{19}c_{23}-2}(\log d_L)^{2-4c_{19}c_{23}}
+ \frac{c_{15}^{\prime}}{d_L^{2c_{23}-c_{10}}}
+ \frac{2\alpha_4 c_{23}^{\frac{1}{2}}}{\log 3} \frac{(\log d_L)^{\frac{3}{2}}}{d_L^{c_{23}-c_{10}}}.$$
If $c_{23}=179$, then $$\frac{9}{10}c_{23} > \epsilon_4.$$
The result follows.
\end{pf}

Lemma \ref{Lem-PositivityOfSumOfPCK1} and \ref{Lem-PositivityOfSumOfPCK2} yield Theorem \ref{Thm-WithoutGRH}.
\\
\\
\textbf{Acknowledgements.~~~} The authors would like to thank the
referee for useful suggestions to improve Section
\ref{Sec-DHPhenomenon}. The first author was supported by Basic
Science Research Program through the National Research Foundation
of Korea funded by the Ministry of Education(NRF-2013R1A1A2061231)
and a Korea University Grant. The second author was supported by
NRF-2013R1A1A2007418.

\end{document}